\documentclass[11pt]{article} % Stat tech report
\newcommand{\doctype}{TECH}

\usepackage{ifthen}

\usepackage{amsfonts,amsmath,amssymb,amsthm}
\usepackage{amsfonts,graphicx,subfigure}
\usepackage{amsfonts,amsmath,amssymb}
\usepackage{epic,eepic,eepicemu}
\usepackage{epsf}
\usepackage{epsfig}
\usepackage{fancyhdr}
\usepackage{graphics}
\usepackage{psfrag,latexsym}

\RequirePackage{natbib}

\ifthenelse{\equal{\doctype}{AOS}}{
\RequirePackage{hypernat}
}
{
}

\ifthenelse{\equal{\doctype}{TECH}}
{
\setlength{\topmargin}{0 in}
\setlength{\textwidth}{6 in}
\setlength{\textheight}{8.5 in}
\setlength{\oddsidemargin}{0.5 in}
\usepackage{natbib}
\bibliographystyle{plainnat}
\bibpunct{(}{)}{,}{a}{,}{,}
}{}

\ifthenelse{\equal{\doctype}{AOS}}
{
\bibliographystyle{plainnat}
\startlocaldefs
}{}
\newcommand{\neglog}{\ensuremath{\ell}}

\newcommand{\Nset}{\ensuremath{S}}
\newcommand{\Nsetcomp}{{\ensuremath{{\Nset^{c}}}}}

\newcommand{\assdep}{A1}
\newcommand{\assinc}{A2}
\newcommand{\assgro}{A3}

\newcommand{\beq}{\begin{eqnarray*}}
\newcommand{\eeq}{\end{eqnarray*}}
\newcommand{\beqn}{\begin{eqnarray}}
\newcommand{\eeqn}{\end{eqnarray}}

\newcommand{\dualvec}{\ensuremath{\estim{z}}}
\newcommand{\obsnoise}{\ensuremath{W^\numobs}}
\newcommand{\Rem}{\ensuremath{R}}
\newcommand{\Qobs}{\ensuremath{\Qmat^{\numobs}}}
\newcommand{\Qstar}{\ensuremath{\Qmat^{*}}}

\let\hat\widehat

\def\x{x}
\def\xi{x^{(i)}}

\def\N{\mathcal{N}}

\newcommand{\degmax}{\ensuremath{\nodedeg}}
\newcommand{\nodedeg}{\ensuremath{d}}
\newcommand{\myparagraph}[1]{\noindent \paragraph{#1}}
\newcommand{\thetastar}{\ensuremath{\theta^*}}
\newcommand{\pdim}{\mdim}
\newcommand{\myeigmin}{\ensuremath{\Lambda_{min}}}
\newcommand{\myeigmax}{\ensuremath{\Lambda_{max}}}

\newcommand{\Qmat}{\ensuremath{Q}}
\newcommand{\Sset}{\ensuremath{S}}

\newcommand{\Sbar}{\ensuremath{{S^{c}}}}
\newcommand{\sgn}{\ensuremath{\operatorname{sgn}}}

\newcommand{\Cmin}{\ensuremath{C_{min}}}

\newcommand{\numobs}{\ensuremath{n}}

\newcommand{\mdim}{\ensuremath{p}}

\theoremstyle{plain}

\newtheorem{theo}{Theorem}[section]

\newtheorem{lem}{Lemma}[section]
\newtheorem{prop}{Proposition}[section]
\newtheorem{cor}{Corollary}[section]

\theoremstyle{definition} 

\newtheorem{nota}{Notation}[section]
\newtheorem{de}{Definition}[section]
\newtheorem{exa}{Example}[section]
\newtheorem{as}{Assumption}[section]
\newtheorem{alg}{Algorithm}[section]

\newcommand{\btheo}{\begin{theo}}
\newcommand{\bde}{\begin{de}}
\newcommand{\ble}{\begin{lem}}
\newcommand{\bpr}{\begin{prop}}
\newcommand{\bno}{\begin{nota}}
\newcommand{\bex}{\begin{exa}}
\newcommand{\bcor}{\begin{cor}}
\newcommand{\spro}{\begin{proof}}
\newcommand{\bas}{\begin{as}}
\newcommand{\balg}{\begin{alg}}

\newcommand{\etheo}{\end{theo}}
\newcommand{\ede}{\end{de}}
\newcommand{\ele}{\end{lem}}
\newcommand{\epr}{\end{prop}}
\newcommand{\eno}{\end{nota}}
\newcommand{\eex}{\end{exa}}
\newcommand{\ecor}{\end{cor}}
\newcommand{\fpro}{\end{proof}}
\newcommand{\eas}{\end{as}}
\newcommand{\ealg}{\end{alg}}

\theoremstyle{plain}

\newtheorem{theos}{Theorem}
\newtheorem{props}{Proposition}
\newtheorem{lems}{Lemma}
\newtheorem{cors}{Corollary}

\theoremstyle{definition}
\newtheorem{exas}{Example}
\newtheorem{algs}{Algorithm}
\newtheorem{asss}{Assumption}
\newtheorem{defns}{Definition}

\newcommand{\btheos}{\begin{theos}}
\newcommand{\etheos}{\end{theos}}
\newcommand{\bprops}{\begin{props}}
\newcommand{\eprops}{\end{props}}
\newcommand{\bdes}{\begin{defns}}
\newcommand{\edes}{\end{defns}}
\newcommand{\blems}{\begin{lems}}
\newcommand{\elems}{\end{lems}}
\newcommand{\bcors}{\begin{cors}}
\newcommand{\ecors}{\end{cors}}
\newcommand{\bexs}{\begin{exas}}
\newcommand{\eexs}{\end{exas}}
\newcommand{\balgs}{\begin{algs}}
\newcommand{\ealgs}{\end{algs}}
\newcommand{\bass}{\begin{asss}}
\newcommand{\eass}{\end{asss}}

\newenvironment{carlist}
 {\begin{list}{$\bullet$}
 {\setlength{\topsep}{0in} \setlength{\partopsep}{0in}
  \setlength{\parsep}{0in} \setlength{\itemsep}{\parskip}
  \setlength{\leftmargin}{0.07in} \setlength{\rightmargin}{0.08in}
  \setlength{\listparindent}{0in} \setlength{\labelwidth}{0.08in}
  \setlength{\labelsep}{0.1in} \setlength{\itemindent}{0in}}}
 {\end{list}}

\newcommand{\bcar}{\begin{carlist}}
\newcommand{\ecar}{\end{carlist}}

\newenvironment{carliste}
 {\begin{list}x
 {\setlength{\topsep}{0in} \setlength{\partopsep}{0in}
  \setlength{\parsep}{0in} \setlength{\itemsep}{\parskip}
  \setlength{\leftmargin}{0.07in} \setlength{\rightmargin}{0.08in}
  \setlength{\listparindent}{0in} \setlength{\labelwidth}{0.08in}
  \setlength{\labelsep}{0.1in} \setlength{\itemindent}{0in}}}
 {\end{list}}

\newcommand{\bcare}{\begin{carliste}}
\newcommand{\ecare}{\end{carliste}}

\newcommand{\Prob}{\ensuremath{\mathbb{P}}}

\makeatletter
\long\def\@makecaption#1#2{
        \vskip 0.8ex
        \setbox\@tempboxa\hbox{\small {\bf #1:} #2}
        \parindent 1.5em  %% How can we use the global value of this???
        \dimen0=\hsize
        \advance\dimen0 by -3em
        \ifdim \wd\@tempboxa >\dimen0
                \hbox to \hsize{
                        \parindent 0em
                        \hfil 
                        \parbox{\dimen0}{\def\baselinestretch{0.96}\small
                                {\bf #1.} #2
                                } 
                        \hfil}
        \else \hbox to \hsize{\hfil \box\@tempboxa \hfil}
        \fi
        }
\makeatother

\long\def\comment#1{}

\makeatletter
\def\@cite#1#2{[\if@tempswa #2 \fi #1]}
\makeatother

\makeatletter
\long\def\barenote#1{
    \insert\footins{\footnotesize
    \interlinepenalty\interfootnotelinepenalty 
    \splittopskip\footnotesep
    \splitmaxdepth \dp\strutbox \floatingpenalty \@MM
    \hsize\columnwidth \@parboxrestore
    {\rule{\z@}{\footnotesep}\ignorespaces
      #1\strut}}}
\makeatother
\newcommand{\bit}{\begin{itemize}}
\newcommand{\eit}{\end{itemize}}
\newcommand{\ben}{\begin{enumerate}}
\newcommand{\een}{\end{enumerate}}

\newcommand{\bear}{\begin{eqnarray}}
\newcommand{\eear}{\end{eqnarray}}

\newcommand{\df}{\ensuremath{d}}

\newcommand{\eparams}{{\ensuremath{{\eparam^{*}}}}}

\newcommand{\eparambar}{\ensuremath{\bar{\eparam}}}

\newcommand{\eparam}{\ensuremath{\theta}}

\newcommand{\statesp}{\ensuremath{{\scr{X}}}}

\newcommand{\order}{{\mathcal{O}}}

\newcommand{\graph}{\ensuremath{G}}

\newcommand{\vertex}{\ensuremath{V}}
\newcommand{\edge}{\ensuremath{E}}

\newcommand{\bk}{\ensuremath{\backslash}}

\newcommand{\Exs}{\ensuremath{{\mathbb{E}}}}

\newcommand{\estart}{\begin{equation}}
\newcommand{\eend}{\end{equation}}

\newcommand{\widgraph}[2]{\includegraphics[keepaspectratio,width=#1]{#2}}

\newcommand{\scr}[1]{\ensuremath{\mathcal{#1}}}

\newcommand{\defn}{\ensuremath{:  =}}

\newcommand{\bec}{\begin{center}}
\newcommand{\enc}{\end{center}}

\newcommand{\beit}{\begin{itemize}}
\newcommand{\enit}{\end{itemize}}

\newcommand{\been}{\begin{enumerate}}
\newcommand{\enen}{\end{enumerate}}

\newcommand{\comsl}{\begin{slide}}
\newcommand{\comspor}{\begin{slide*}}

\newcommand{\comsld}[2]{\begin{slide}[#1,#2]}
\newcommand{\comspord}[2]{\begin{slide*}[#1,#2]}

\newcommand{\mendsl}{\end{slide}}
\newcommand{\mendspo}{\end{slide*}}

\newcommand{\estim}[1]{\ensuremath{\widehat{#1}}}

\newcommand{\wtil}[1]{\ensuremath{\widetilde{#1}}}

\newcommand{\real}{\ensuremath{{\mathbb{R}}}}

\DeclareMathOperator{\sign}{sign}

\newcommand{\Wnoise}{\ensuremath{W}}
\newcommand{\mythetahat}{\widehat{\theta}}

\newcommand{\mprob}{\ensuremath{\mathbb{P}}}

\newcommand{\Rvar}{\ensuremath{X}}

\newcommand{\eparamhat}{\estim{\eparam}}

\newcommand{\varfun}[2]{\ensuremath{\eta(#1; #2)}}
\newcommand{\Dmax}{\ensuremath{D_{\operatorname{max}}}}

\newcommand{\Data}{\ensuremath{\{\xi\}}}

\newcommand{\matnorm}[3]{|\!|\!| #1 | \! | \!|_{{#2}}}
\newcommand{\matsnorm}[2]{|\!|\!| #1 | \! | \!|_{{#2}}}

\newcommand{\mutinc}{\ensuremath{\alpha_{\operatorname{{\scriptstyle{}}}}}}

\newcommand{\Zi}{\ensuremath{Z^{(i)}}}
\newcommand{\myZ}{\ensuremath{Z}}

\newcommand{\Term}{\ensuremath{T}}

\newcommand{\unicon}{\ensuremath{K}}

\newcommand{\regpar}{\ensuremath{\lambda}}

\newcommand{\thetamin}{\ensuremath{\theta^*_{\operatorname{min}}}}

\newcommand{\vnots}[1]{\ensuremath{\backslash #1}}

\newcommand{\myxsam}[1]{\ensuremath{x^{(#1)}}}

\newcommand{\muhat}{\ensuremath{\estim{\mu}}}
\newcommand{\Logfun}[2]{\ensuremath{f(#1; #2)}}

\newcommand{\hatthetas}[1]{\ensuremath{\estim{\theta}^\numobs_{#1}}}

\newcommand{\hatNsign}{\estim{\N}_\pm}

\newcommand{\svert}{\ensuremath{r}}

\newcommand{\zhat}{\ensuremath{\estim{z}}}

\newcommand{\signvec}{\ensuremath{{\vec{v}}}}
\newcommand{\psignvec}{v}

\newcommand{\GoodEvent}{\ensuremath{\mathcal{M}}}

\newcommand{\thetahat}{\ensuremath{\estim{\eparam}}}

\newcommand{\contpar}{\ensuremath{\beta}}

\newcommand{\Nsign}{\ensuremath{\N_{\pm}}}

\newcommand{\zsam}[1]{\ensuremath{Z^{(#1)}}}

\newcommand{\pert}{\ensuremath{u}}
\newcommand{\mydefn}{\ensuremath{: \, =}}
\newcommand{\Mrad}{\ensuremath{M}}

\def\beq{\begin{eqnarray*}}
\def\eeq{\end{eqnarray*}}

\def\tparam{\eparam^*}

\def\Hess{\nabla^2}

\ifthenelse{\equal{\doctype}{AOS}}
{
\endlocaldefs
}{}

\begin{document}

\ifthenelse{\equal{\doctype}{TECH}}
{

\begin{center}

{\bf{\Large{High-Dimensional Graphical Model Selection Using
$\ell_1$-Regularized Logistic Regression}}}

\vspace*{.2in}

\begin{tabular}{ccc}
Pradeep Ravikumar$^\dagger$ & \hspace*{.2in} & Martin J. Wainwright$^{\ast, \dagger}$ \\
\texttt{pradeepr@stat.berkeley.edu} & &
\texttt{wainwrig@stat.berkeley.edu} 
\end{tabular}

\vspace*{.1in}

\begin{tabular}{c}
Department of Statistics$^{\dagger}$, and Department of EECS$^\ast$ \\
UC Berkeley, Berkeley, CA 94720
\end{tabular}

\vspace*{.2in}

\begin{tabular}{c}
John D. Lafferty \\
\texttt{lafferty@cs.cmu.edu}
\end{tabular}

\vspace*{.1in}

\begin{tabular}{c}
Departments of Computer Science and Machine Learning \\
Carnegie Mellon University  \\
Pittsburgh, PA 15213 
\end{tabular}

\vspace*{.1in}

\begin{tabular}{c}
Technical Report, Department of Statistics, UC Berkekley \\
April 19, 2008
\end{tabular}

\vspace*{.1in}

\begin{abstract}
We consider the problem of estimating the graph structure associated
with a discrete Markov random field.  We describe a method based on
$\ell_1$-regularized logistic regression, in which the neighborhood of
any given node is estimated by performing logistic regression subject
to an $\ell_1$-constraint.  Our framework applies to the
high-dimensional setting, in which both the number of nodes $p$ and
maximum neighborhood sizes $d$ are allowed to grow as a function of
the number of observations $n$.  Our main results provide sufficient
conditions on the triple $(n, p, d)$ for the method to succeed in
consistently estimating the neighborhood of every node in the graph
simultaneously.  Under certain assumptions on the population Fisher
information matrix, we prove that consistent neighborhood selection
can be obtained for sample sizes $n = \Omega(d^3 \log p)$, with the
error decaying as $\order(\exp(-C n/d^3))$ for some constant $C$.  If
these same assumptions are imposed directly on the sample matrices, we
show that $n = \Omega(d^2 \log p)$ samples are sufficient.
\end{abstract}

\vspace*{.1in}
\end{center}


}{}

\ifthenelse{\equal{\doctype}{AOS}}
{
\begin{frontmatter}

\title{High-Dimensional Graphical Model Selection Using
$\ell_1$-Regularized Logistic Regression}

\runtitle{High-Dimensional Graphical Model Selection}

\author{\fnms{Pradeep} \snm{Ravikumar}\ead[label=e1]{pradeepr@stat.berkeley.edu}\thanksref{t1,t2,t3}},
\author{\fnms{Martin J.} \snm{Wainwright} \ead[label=e2]{wainwrig@stat.berkeley.edu}\thanksref{t2}}
\affiliation{University of California, Berkeley}

\author{\fnms{John D.} \snm{Lafferty} \ead[label=e3]{lafferty@cs.cmu.edu}\thanksref{t1}}
\affiliation{Carnegie Mellon University}
\address{Berkeley, CA 94720 \printead{e1,e2}\\
Pittsburgh, PA 15213 \printead{e3}}

%PR and JL
\thankstext{t1}{Supported in part by NSF grants IIS-0427206 and CCF-0625879}
%PR and MJW
\thankstext{t2}{Supported in part by NSF grants DMS-0605165 and CCF-0545862}
%PR
\thankstext{t3}{Supported in part by a Siebel Scholarship}

\runauthor{Ravikumar, Wainwright and Lafferty}

\begin{abstract}
We consider the problem of estimating the graph structure associated
with a discrete Markov random field.  We describe a method based on
$\ell_1$-regularized logistic regression, in which the neighborhood of
any given node is estimated by performing logistic regression subject
to an $\ell_1$-constraint.  Our framework applies to the
high-dimensional setting, in which both the number of nodes $p$ and
maximum neighborhood sizes $d$ are allowed to grow as a function of
the number of observations $n$.  Our main results provide sufficient
conditions on the triple $(n, p, d)$ for the method to succeed in
consistently estimating the neighborhood of every node in the graph
simultaneously.  Under certain assumptions on the population Fisher
information matrix, we prove that consistent neighborhood selection
can be obtained for sample sizes $n = \Omega(d^3 \log p)$, with the
error decaying as $\order(\exp(-C n/d^3))$ for some constant $C$.  If
these same assumptions are imposed directly on the sample matrices, we
show that $n = \Omega(d^2 \log p)$ samples are sufficient.
\end{abstract}



\end{frontmatter}

}{}

\ifthenelse{\equal{\doctype}{JMLR}}
{
\title{High-Dimensional Graphical Model Selection Using
$\ell_1$-Regularized Logistic Regression}

\author{\name Pradeep Ravikumar \email pradeepr@stat.berkeley.edu\\
\addr Department of Statistics\\ University of California\\ Berkeley,
CA 94720-1776, USA \AND \name Martin J. Wainwright \email
wainwrig@stat.berkeley.edu\\ \addr Department of Statistics and
Department of EECS\\ University of California\\ Berkeley, CA
94720-1776, USA \AND \name John D. Lafferty \email
lafferty@cs.cmu.edu\\ \addr Computer Science Department and Machine Learning Department
\\ Carnegie Mellon University\\ Pittsburgh, PA
15213, USA}

\maketitle

\begin{abstract}
We consider the problem of estimating the graph structure associated
with a discrete Markov random field.  We describe a method based on
$\ell_1$-regularized logistic regression, in which the neighborhood of
any given node is estimated by performing logistic regression subject
to an $\ell_1$-constraint.  Our framework applies to the
high-dimensional setting, in which both the number of nodes $p$ and
maximum neighborhood sizes $d$ are allowed to grow as a function of
the number of observations $n$.  Our main results provide sufficient
conditions on the triple $(n, p, d)$ for the method to succeed in
consistently estimating the neighborhood of every node in the graph
simultaneously.  Under certain assumptions on the population Fisher
information matrix, we prove that consistent neighborhood selection
can be obtained for sample sizes $n = \Omega(d^3 \log p)$, with the
error decaying as $\order(\exp(-C n/d^3))$ for some constant $C$.  If
these same assumptions are imposed directly on the sample matrices, we
show that $n = \Omega(d^2 \log p)$ samples are sufficient.
\end{abstract}

\vskip5pt

}{}

%%%%%%%%%%%%%%%%%%%%%%%%%%%%%%%%%%%%%%%%%%%%%%%%%%%%%%%%%%%%%%%%%%%%%%%%%%%%%%%%%%

\ifthenelse{\equal{\doctype}{TECH}}
{
%\begin{onehalfspace}
%\begin{doublespace}
\renewcommand{\baselinestretch}{1}
% Hack to make renewcommand take effect immediately (if not in preamble,
% only takes effect after a switch in font size)
\small \normalsize 
}
%%%%%%%
{
%\typeout{NO FORMATTING CHANGE FOR JMLR OR AOS}
}

%%%%%%%%%%%%%%%%%%%%%%%%%%%%%%%%%%%%%%%%%%%%%%%%%%%%%%%%%%%%%%%%%%%%%%%%%%%%%
\section{Introduction}

Undirected graphical models, also known as Markov random fields
(MRFs), are used in a variety of domains, including artificial
intelligence, natural language processing, image analysis, statistical
physics, and spatial statistics, among others.  A Markov random field
(MRF) is specified by an undirected graph $\graph = (\vertex, \edge)$,
with vertex set $\vertex = \{1, 2, \ldots, \pdim \}$ and edge set
$\edge \subset \vertex \times \vertex$.  The structure of this graph
encodes certain conditional independence assumptions among subsets of
the $\pdim$-dimensional discrete random variable $\Rvar = (\Rvar_1,
\Rvar_2,\ldots, \Rvar_\mdim)$, where variable $\Rvar_i$ is associated
with vertex $i \in \vertex$.  A fundamental problem is the
\emph{graphical model selection problem}: given a set of $\numobs$
samples $\{ \myxsam{1}, \myxsam{2}, \ldots, \myxsam{\numobs} \}$ from
a Markov random field, estimate the structure of the underlying graph.
The \emph{sample
complexity} of such an estimator is the minimal number of samples $\numobs$, as a function
of the graph size $\pdim$ and possibly other parameters such as the
maximum node degree $\degmax$, required for the probability of correct
identification of the graph to converge to one.  Another 
important property of any model selection procedure is its
\emph{computational complexity}.

Due to both its importance and difficulty, structure learning in
random fields has attracted considerable attention.  The absence of an edge
in a graphical model encodes a conditional independence assumption.
Constraint-based approaches~\citep{spirtes:00} estimate these
conditional independencies from the data using hypothesis testing, and
then determine a graph that most closely represents those
independencies. Each graph represents a model class of graphical
models; learning a graph then is a model class selection
problem. Score-based approaches combine a metric for the complexity of
the graph, with a goodness of fit measure of the graph to the data
(for instance, log-likelihood of the maximum likelihood parameters
given the graph), to obtain a \emph{score} for each graph. The score
is then used together with a search procedure that generates candidate
graph structures to be scored. The number of graph structures grows
super-exponentially, however, and~\cite{chickering:95} shows that this
problem is in general NP-hard. 

A complication for undirected graphical models is that typical score
metrics involve the normalization constant (also called the partition
function) associated with the Markov random field, which is
intractable (\#P) to compute for general undirected models.  The space
of candidate structures in scoring based approaches is thus typically
restricted to either directed models---Bayesian networks---or simple
undirected graph classes such as trees~\citep{chowliu:68},
polytrees~\citep{dasgupta:99} and hypertrees~\citep{srebro:03}.
\cite{AbbKolNg06} propose a method for learning factor graphs based on
local conditional entropies and thresholding, and analyze its behavior
in terms of Kullback-Leibler divergence between the fitted and true
models.  They obtain a sample complexity that grows logarithmically in
the graph size $\pdim$, but the computational complexity grows at
least as quickly as $\order(\pdim^{\degmax+1})$, where $\degmax$ is
the maximum neighborhood size in the graphical model. This order of
complexity arises from the fact that for each node, there are
$\binom{\pdim}{\degmax} = \order(\pdim^\degmax)$ possible
neighborhoods of size $\degmax$ for a graph with $\pdim$ vertices.
\cite{Csiszar:06} show consistency of a method that uses
pseudo-likelihood and a modification of the BIC criterion, but this
also involves a prohibitively expensive search.

In work subsequent to the initial conference version of this work
\citep{WaiRavLaf06}, other researchers have also studied the problem
of model selection in discrete Markov random fields.  For the special
case of bounded degree models, \cite{Bresler08} describe a simple
search-based method, and prove under relatively mild assumptions that
it can recover the graph structure with $\Theta(\degmax \log \pdim)$
samples.  However, in the absence of additional restrictions, the
computational complexity of the method is $\order(\pdim^{\degmax+1})$.
\noindent\cite{SanWai08} analyze the information-theoretic limits of
graphical model selection, providing both upper and lower bounds on
various model selection procedures, but these methods also have
prohibitive computational cost.

The main contribution of this paper is to present and analyze the
computational and sample complexity of a simple method for graphical
model selection.  Our analysis is high-dimensional in nature, meaning
that both the model dimension $\pdim$ as well as the maximum
neighborhood size $\degmax$ may tend to infinity as a function of the
size $\numobs$.  Our main result shows that under mild assumptions,
consistent neighborhood selection is possible with sample complexity
$\numobs = \Omega(\degmax^3 \log \pdim)$ and computational complexity
$\order(\max\{\numobs, \pdim \} \pdim^3 )$, when applied to any graph
with $\pdim$ vertices and maximum degree $\degmax$.  The basic
approach is straightforward: it involves performing
$\ell_1$-regularized logistic regression of each variable on the
remaining variables, and then using the sparsity pattern of the
regression vector to infer the underlying neighborhood structure.  

The technique of $\ell_1$ regularization for estimation of sparse
models or signals has a long history in many fields; for instance,
see~\cite{Tropp:06} for one survey.  A surge of recent work has shown
that $\ell_1$-regularization can lead to practical algorithms with
strong theoretical guarantees
(e.g.,~\cite{CandesTao06,Donoho:Elad:03,Meinshausen:06,ng:04,Tropp:06,Wainwright06a_aller,ZhaoYu06}).
Despite the well-known computational intractability of discrete MRFs,
our method is computationally efficient; it involves neither computing
the normalization constant (or partition function) associated with the
Markov random field, nor a combinatorial search through the space of
graph structures. Rather, it requires only the solution of standard
convex programs, with an overall computational complexity of order
$\order(\max \{\pdim, \numobs\} \, \pdim^3)$ \citep{KohKimBoy07}, and
is thus well-suited to high dimensional problems.  Conceptually, like
the work of~\cite{Meinshausen:06} on covariance selection in Gaussian
graphical models, our approach can be understood as using a type of
pseudo-likelihood, based on the local conditional likelihood at each
node.  In contrast to the Gaussian case, where the exact maximum
likelihood estimate can be computed exactly in polynomial time, this
use of a surrogate loss function is essential for discrete Markov
random fields, given the intractability of computing the exact
likelihood.

The remainder of this paper is organized as follows.  We begin in
Section~\ref{SecBackground} with background on discrete graphical
models, the model selection problem, and logistic regression.  In
Section~\ref{SecOutline}, we state our main result, develop some of
its consequences, and provide a high-level outline of the proof.
Section~\ref{SecFixedDesign} is devoted to proving a result under
stronger assumptions on the sample Fisher information matrix itself,
whereas Section~\ref{SecPopulation} provides concentration results
linking the population matrices to the sample versions.  In
Section~\ref{SecExperiments}, we provide some experimental results to
illustrate the practical performance of our method, and the close
agreement between theory and practice, and we conclude in
Section~\ref{SecDiscussion}.

\myparagraph{Notation:} For the convenience of the reader, we
summarize here notation to be used throughout the paper.  We use the
following standard notation for asymptotics: we write $f(n) =
\order(g(n))$ if $f(n) \leq \unicon g(n)$ for some constant $\unicon <
\infty$, and $f(n) = \Omega(g(n))$ if $f(n) \geq \unicon' g(n)$ for
some constant $\unicon' > 0$.  The notation $f(n) = \Theta(g(n))$
means that $f(n) = \order(g(n))$ and $f(n) = \Omega(g(n))$. Given a
vector $x \in \real^\df$ and parameter $q \in [1, \infty]$, we use
$\|\x\|_q$ to denote the usual $\ell_q$ norm.  Given a matrix $X \in
\real^{a \times b}$ and parameter $q \in [1, \infty]$, we use
$\matsnorm{X}{q}$ to denote the induced matrix-operator norm (viewed
as a mapping from $\ell_q^b \rightarrow \ell_q^a$); see \cite{Horn85}.
Two examples of particular importance in this paper are the spectral
norm $\matsnorm{X}{2}$, corresponding to the maximal singular value of
$X$, and the $\ell_\infty$ matrix norm, given by $\matsnorm{X}{\infty}
=\max \limits_{j=1, \ldots, a} \sum_{k=1}^b |X_{jk}|$.  We make use of
the bound $\matsnorm{X}{\infty} \leq \sqrt{a} \matsnorm{X}{2}$, for
any square matrix $X \in \real^{a \times a}$.

\section{Background and problem formulation}
\label{SecBackground}

We begin by providing some background on Markov random fields,
defining the problem of graphical model selection, and describing our
method based on neighborhood logistic regression.

\subsection{Markov random fields}

Given an undirected graph $\graph$ with vertex set $\vertex = \{1,
\ldots, \pdim \}$ and edge set $\edge$, a Markov random field (MRF)
consists of random vector $X = (X_1, X_2, \ldots, X_\pdim)$, where the
random variable $X_s$ is associated with vertex $s \in \vertex$.  The
random vector $X \in \statesp^\pdim$ is said to be pairwise Markov
with respect to the graph if its probability distribution factorizes
as $\mprob(x) \propto \exp \left \{ \sum_{(s,t) \in \edge}
\phi_{st}(x_s, x_t) \right\}$, where each $\phi_{st}$ is a mapping
from pairs $(x_s, x_t) \in \statesp_s \times \statesp_t$ to the real
line.  An important special case is the Ising model, in which $X_s \in
\{-1,1\}$ for each vertex $s \in \vertex$, and $\phi_{st}(x_s, x_t) =
\eparam_{st} x_s x_t$ for some parameter $\eparam^*_{st} \in \real$, so
that the distribution takes the form
\begin{equation}
\label{EqnIsing}
\mprob_{\eparams}(x) = \frac{1}{Z(\eparams)} \exp \left \{\sum_{(s,t)\in
\edge} \eparam^*_{st} x_s x_t \right \}.
\end{equation}
The partition function $Z(\eparams)$ ensures that the distribution
sums to one.  The Ising model has proven useful in many domains,
including statistical physics, where it describes the behavior of
gases or magnets, in computer vision for image segmentation, and in
social network analysis.

\subsection{Graphical model selection}

Suppose that we are given a collection $\Data = \{\myxsam{1}, \ldots,
\myxsam{\numobs} \}$ of $\numobs$ samples, where each $\pdim$-dimensional
vector $\myxsam{i}$ is drawn in an i.i.d. manner from a distribution
$\mprob_{\eparams}$ of the form~\eqref{EqnIsing}.  It is convenient to
view the parameter vector $\eparams$ as a
$\binom{\pdim}{2}$-dimensional 
vector, indexed by pairs of distinct vertices, but
non-zero \emph{if and only if} the vertex pair $(s,t)$ belongs to the
unknown edge set $\edge$ of the underlying graph $\graph$.  The goal
of \emph{graphical model selection} is to infer the edge set $\edge$
of the graphical model defining the probability distribution
that generates the samples.  In this paper, we study the slightly stronger
criterion of \emph{signed edge recovery}: in particular, given a
graphical model with parameter $\eparams$, we define the edge sign
vector
\begin{eqnarray}
\label{EqnSignedEdge}
\edge^* & \defn & \begin{cases} \sign(\eparam^*_{st}) & \mbox{if $(s,t) \in
\edge$} \\ 0 & \mbox{otherwise.}
				     \end{cases}
\end{eqnarray}
Note that the weaker graphical model selection problem amounts to
recovering the vector $|\edge^*|$ of absolute values.

The classical notion of statistical consistency applies to the
limiting behavior of an estimation procedure as the sample size
$\numobs$ goes to infinity, with the model size $\pdim$ itself
remaining fixed.  In many contemporary applications of graphical
models (e.g., gene microarrays, social networks etc.), the model
dimension $\pdim$ is comparable or larger than the sample size
$\numobs$, so that the relevance of such ``fixed $\pdim$'' asymptotics
is doubtful.  Accordingly, the goal of this paper is to develop
the broader notion of \emph{high-dimensional consistency}, in which
both the model dimension and the sample size are allowed to increase,
and we study the scaling conditions under which consistent model selection
is achievable.

More precisely, we consider sequences of graphical model selection
problems, indexed by the sample size $\numobs$, number of vertices
$\pdim$, and maximum node degree $\degmax$.  We assume that the sample
size $\numobs$ goes to infinity, and both the problem dimension $\pdim
= \pdim(\numobs)$ and $\degmax = \degmax(\numobs)$ may also scale as a
function of $\numobs$.  The setting of fixed $\pdim$ or
$\degmax$ is covered as a special case. Let
$\estim{\edge}_\numobs$ be an estimator of the signed edge pattern
$\edge^*$, based on the $\numobs$ samples.  Our goal is to establish
sufficient conditions on the scaling of the triple $(\numobs, \pdim,
\degmax)$ such that our proposed estimator is consistent in the sense
that
\begin{eqnarray*}
\mprob \left[\widehat{\edge}_\numobs = \edge^* \right] & \rightarrow &
1 \qquad \mbox{as $\numobs \rightarrow +\infty$}.
\end{eqnarray*}
We sometimes call this property \textit{sparsistency}, as a shorthand for
consistency of the sparsity pattern of the parameter $\theta^*$.

\subsection{Neighborhood-based logistic regression}

Note that recovering the signed edge vector $\edge^*$ of an undirected
graph $\graph$ is equivalent to recovering, for each vertex $\svert
\in \vertex$, its \emph{neighborhood set} $\N(\svert) \defn \{t \in
\vertex \, \mid \, (\svert,t) \in \edge\}$, along with the correct
signs $\sign(\eparam^*_{\svert t})$ for all $t \in \N(\svert)$.  To
capture both the neighborhood structure and sign pattern, we define
the \emph{signed neighborhood set}
\begin{eqnarray}
\label{EqnSignedNeigh}
\Nsign(\svert) & \defn & \left \{ \sign(\eparam^*_{\svert t}) \, t \;
\mid \; t \in \N(s) \right \}.
\end{eqnarray}
The next step is to observe that this signed neighborhood set
$\Nsign(\svert)$ can be recovered from the sign-sparsity pattern of
the $(\pdim-1)$-dimensional subvector of parameters
\begin{eqnarray*}
\eparam^*_{\bk \svert} & \defn & \left \{ \eparam^*_{\svert u}, \; u \in
\vertex \bk \svert \right \}
\end{eqnarray*}
associated with vertex $\svert$.  In order to estimate this vector
$\eparam^*_{\bk \svert}$, we consider the structure of the conditional
distribution of $X_\svert$ given the other variables
$X_{\vnots{\svert}} = \{X_t \, \mid \, t \in \vertex \bk \{\svert\}
\}$.  A simple calculation shows that under the
model~\eqref{EqnIsing}, this conditional distribution takes the form
\begin{eqnarray}
\label{EqnCondDist}
\mprob_{\eparams}(x_\svert \, \mid x_{\vnots{\svert}}) & = & \frac{
 \exp \left( 2 x_\svert \left[\sum_{t \in \vertex \bk \svert}
 \eparam^*_{\svert t} x_t \right] \right)}{ \exp \left(2 x_\svert[
 \sum_{t \in \vertex \bk \svert} \eparam^*_{\svert t} x_t] \right) + 1
 }.
\end{eqnarray}
Thus, the variable $X_\svert$ can be viewed as the response variable
in a logistic regression in which all of the other variables
$X_{\vnots{\svert}}$ play the role of the covariates.

With this set-up, our method for estimating the sign-sparsity pattern
of the regression vector $\eparam^*_{\bk \svert}$ (and hence the
neighborhood structure $\Nsign(\svert)$) is based on computing the
$\ell_1$-regularized logistic regression of $X_\svert$ on the other
variables $X_{\vnots{\svert}}$.  Explicitly, given a set of $\numobs$
i.i.d. samples $\{\myxsam{1}, \myxsam{2}, \ldots, \myxsam{\numobs}
\}$, this regularized regression problem is a convex program, of the
form 
$$\min \limits_{\eparam_{\vnots{\svert}} \in \real^{\pdim -1} }
\left \{ \neglog(\eparam; \Data) + \regpar_\numobs
\|\eparam_{\vnots{\svert}}\|_1 \right \}$$, where $\regpar_\numobs > 0$
is a regularization parameter, to be specified by the user, and
\begin{eqnarray}
\neglog(\eparam; \Data) & \defn & -\frac{1}{\numobs}
\sum_{i=1}^\numobs \log \mprob_{\eparam}(\myxsam{i}_\svert \, \mid
\myxsam{i}_{\vnots{\svert}})
\end{eqnarray}
is the rescaled negative log likelihood.  (The rescaling factor
$1/\numobs$ in this definition is for later theoretical convenience.)
Following some algebraic manipulation, the regularized negative log
likelihood can be written as
\begin{equation}
\label{EqnRegLikeTwo}
\min \limits_{\eparam_{\vnots{\svert}} \in \real^{\pdim -1} } \left \{
  \frac{1}{\numobs} \sum_{i=1}^\numobs \Logfun{\eparam}{\myxsam{i}} -
  \sum_{u \in \vertex \bk \svert} \eparam_{\svert u} \muhat_{\svert u}
  + \regpar_\numobs \|\eparam_{\vnots{\svert}}\|_1 \right \},
\end{equation}
where 
\begin{eqnarray}
\label{EqnDefnLogFun}
\Logfun{\eparam}{x} & \defn & \log \left(\exp ( \sum_{t \in
\vertex \bk \svert} \eparam_{\svert t} x_t) + \exp (- [\sum_{t \in
\vertex \bk \svert} \eparam_{\svert t} x_t])\right)
\end{eqnarray}
is a rescaled logistic loss, and $\muhat_{\svert u} \defn
\frac{1}{\numobs} \sum_{i=1}^\numobs \myxsam{i}_\svert \myxsam{i}_u$
are empirical moments.  Note the objective
function~\eqref{EqnRegLikeTwo} is convex but not differentiable, due
to the presence of the $\ell_1$-regularizer.  By Lagrangian duality,
the problem~\eqref{EqnRegLikeTwo} can be re-cast as a constrained
problem over the ball $\|\eparam_{\vnots{\svert}}\|_1 \leq
C(\regpar_\numobs)$.  Consequently, by the Weierstrass theorem, the
minimum over $\eparam_{\vnots{s}}$ is always achieved.

Accordingly, let $\hatthetas{\vnots{\svert}}$ be an element of the
minimizing set of problem~\eqref{EqnRegLikeTwo}.  Although
$\hatthetas{\vnots{\svert}}$ need not be unique in general since the
problem~\eqref{EqnRegLikeTwo} need not be strictly convex, our
analysis shows that in the regime of interest, this minimizer
$\hatthetas{\vnots{\svert}}$ is indeed unique.  We use
$\hatthetas{\vnots{\svert}}$ to estimate the signed neighborhood
$\Nsign(\svert)$ according to
\begin{eqnarray}
\label{EqnNeighEst}
\hatNsign(\svert) & \defn & \left \{ \sign(\estim{\eparam}_{\svert u})
u \; \mid \; u \in \vertex \bk \svert, \; \estim{\eparam}_{su} \neq 0
\right \}.
\end{eqnarray}
We say that the full graph $\graph$ is estimated consistently, written
as the event $\{ \estim{\graph} = \graph(\pdim, \degmax) \}$, if
$\hatNsign(\svert) = \Nsign(\svert)$ for all $\svert \in \vertex$.

\section{Method and theoretical guarantees}
\label{SecOutline}

Our main result concerns conditions on the sample size $\numobs$
relative to the parameters of the graphical model---more specifically,
the number of nodes $\mdim$ and maximum node degree $\degmax$---that
ensure that the collection of signed neighborhood
estimates~\eqref{EqnNeighEst}, one for each node $\svert$ of the
graph, agree with the true neighborhoods, so that the full graph
$\graph(\pdim, \degmax)$ is estimated consistently.  In this section,
we begin by stating the assumptions that underlie our main result,
and then give a precise statement of the main result. We then provide a
high-level overview of the key steps involved in its proof, deferring
detail to later sections.  Our analysis proceeds by first establishing
sufficient conditions for correct signed neighborhood recovery---that
is, $\{ \hatNsign(\svert) = \Nsign(\svert) \}$---for some fixed node
$\svert \in \vertex$.  By showing that this neighborhood consistency
is achieved at exponentially fast rates, we can then use a union bound
over all $\pdim$ nodes of the graph to conclude that consistent graph
selection $\{\estim{\graph} = \graph(\pdim, \degmax) \}$ is also
achieved.

\subsection{Assumptions}

Success of our method requires certain assumptions on the structure of
the logistic regression problem.  These assumptions are stated in
terms of the Hessian of the likelihood function $\Exs \{\log
\mprob_{\eparam}[X_\svert \, \mid \, X_{\vnots{\svert}}]\}$, as
evaluated at the true model parameter $\eparam^*_{\vnots{\svert}} \in
\real^{p-1}$. More specifically, for any fixed node $\svert \in
\vertex$, this Hessian is a $(\pdim-1) \times (\pdim-1)$ matrix of the
form
\begin{eqnarray}
\label{EqnDefnQstar}
\Qstar_\svert & \defn & \Exs_{\eparams} \left \{ \nabla^2 \log
\mprob_{\eparams}[X_\svert \, \mid \, X_{\vnots{\svert}}] \right \}.
\end{eqnarray}
For future reference, we calculate the explicit expression
\begin{eqnarray}
\label{EqnDefnQstarExp}
\Qstar_\svert & = & \Exs_{\eparams} \left[\varfun{X}{\eparams}
X_{\vnots{\svert}} X_{\vnots{\svert}}^T \right]
\end{eqnarray}
where
\begin{eqnarray}
\label{EqnDefnVarfun}
\varfun{u}{\eparam} & \defn & \frac{ 4 \exp \left(2 u_\svert
\left[\sum_{t \in \vertex \bk \svert} \eparam_{\svert t} u_t \right]
\right)}{(\exp \left(2 u_\svert [ \sum_{t \in \vertex \bk \svert }
\eparam_{\svert t} u_t] \right) + 1)^2 }
\end{eqnarray}
is the variance function.  Note that the matrix $\Qstar_\svert$ is the
Fisher information matrix associated with the local conditional
probability distribution.  Intuitively, it serves as the counterpart
for discrete graphical models of the covariance matrix $\Exs[X X^T]$
of Gaussian graphical models, and indeed our assumptions are analogous
to those imposed in previous work on the Lasso for Gaussian linear
regression~\citep{Meinshausen:06,Tropp:06,ZhaoYu06}.

In the following we write simply $\Qstar$ for the matrix $\Qstar_\svert$, where
the reference node $\svert$ should be understood implicitly.
Moreover, we use $\Sset \defn \{ (\svert, t) \; \mid \; t \in
\N(\svert) \}$ to denote the subset of indices associated with edges
of $\svert$, and $\Nsetcomp$ to denote its complement.  We use
$\Qstar_{\Sset \Sset}$ to denote the $\degmax \times \degmax$
sub-matrix of $\Qstar$ indexed by $\Sset$.  With this notation, we
state our assumptions:
\myparagraph{[\assdep] Dependency condition:} The subset of the Fisher
information matrix corresponding to the relevant covariates has
bounded eigenvalues: there exists a constant $\Cmin > 0$ such
that 
\begin{equation}
\label{EqnEigBounds}
 \myeigmin(\Qstar_{\Nset \Nset}) \, \geq \; \Cmin.
\end{equation}
Moreover, we require that
$\myeigmax(\Exs_{\eparams}[X_{\vnots{\svert}} X_{\vnots{\svert}}^T])
\leq \Dmax$.  These conditions ensure that the relevant covariates do
not become overly dependent. (As stated earlier, we have suppressed
notational dependence on $\svert$; thus this condition is assumed
to hold for all $\svert \in V$.)

\myparagraph{[\assinc] Incoherence condition:} Our next assumption
captures the intuition that the large number of irrelevant covariates
(i.e., non-neighbors of node $\svert$) cannot exert an overly strong
effect on the subset of relevant covariates (i.e., neighbors of node
$\svert$).  To formalize this intuition, we require the existence of
an $\mutinc \in (0,1]$ such that
\begin{eqnarray}
\label{EqnDefnMutualInco}
\matnorm{\Qmat^*_{\Nsetcomp \Nset} (\Qmat^*_{\Nset
\Nset})^{-1}}{\infty}{\infty} & \leq & 1 - \mutinc.
\end{eqnarray}

\subsection{Statement of main result}
We are now ready to state our main result on the performance of
neighborhood logistic regression for graphical model selection.
Naturally, the limits of model selection are determined by the minimum
value over the parameters $\eparam^*_{\svert t}$ for pairs $(\svert,
t)$ included in the edge set of the true graph.  Accordingly, we
define the parameter
\begin{eqnarray}
\label{EqnDefnTmin}
\thetamin & = & \min_{(\svert, t) \in \edge}  |\eparam^*_{\svert t}|.
\end{eqnarray}
With this definition, we have the following
\btheos
\label{ThmMain}
Consider a sequence of graphs $\{\graph(\pdim, \degmax)\}$ such that
conditions $\assdep$ and $\assinc$ are satisfied by the population
Fisher information matrices $\Qstar$. If the sample size $\numobs$
satisfies
\begin{eqnarray}
\label{EqnGrowthCondition}
\numobs  >  L \degmax^3 \log(\pdim)
\end{eqnarray}
for some constant $L$, and the minimum value $\thetamin$ decays no
faster than $\order(1/\degmax)$, then for the regularization sequence
$\regpar_\numobs = 2 \sqrt{\frac{\log \pdim}{\numobs}}$, the estimated
graph $\estim{\graph}(\regpar_\numobs)$ obtained by neighborhood
logistic regression satisfies
\begin{eqnarray}
\label{EqnSuccess}
\mprob[\estim{\graph}(\regpar_\numobs) = \graph(\pdim, \degmax)] & = &
\order \left(\exp\left(-\unicon \frac{\numobs}{\degmax^3} - 3 \log(\pdim)
\right)\right) \; \rightarrow \; 0
\end{eqnarray}
for some constant $K$.
\etheos

\myparagraph{Remarks:} For model selection in graphical models, one is
typically interested in node degrees $\degmax$ that remain bounded
(e.g., $\degmax = \order(1)$), or that grow only weakly with graph size
(say $\degmax = o(\mdim$).  In such cases, the growth
condition~\eqref{EqnGrowthCondition} allows the number of observations
to be substantially smaller than the graph size, i.e., the ``large
$\pdim$, small $\numobs$'' regime.  In particular, the graph size
$\mdim$ can grow exponentially with the number of observations (i.e,
$\mdim(\numobs) = \exp(\numobs^\alpha)$ for some $\alpha \in (0,1)$.

In terms of the choice of regularization, the sequence
$\regpar_\numobs$ needs to satisfy the following conditions:
\begin{equation*}
\numobs \regpar_\numobs^2 > 2 \log(\pdim), \quad \mbox{and} \quad
\sqrt{\degmax} \regpar_\numobs = \order(\thetamin).
\end{equation*}
Under the growth condition~\eqref{EqnGrowthCondition}, the choice
$\regpar_\numobs = 2\sqrt{\frac{\log \pdim}{\numobs}}$ suffices as
long as $\thetamin$ decays no faster than $\order(1/\degmax)$.

The analysis required to prove Theorem~\ref{ThmMain} can be divided
naturally into two parts.  First, in Section~\ref{SecFixedDesign}, we
prove a result (stated as Proposition~\ref{PropFixed}) for ``fixed
design'' matrices.  More precisely, we show that if the dependence
(\assdep) mutual incoherence (\assinc) conditions hold for the
\emph{sample Fisher information matrix}
\begin{eqnarray}
\label{EqnDefnSampleFisher}
\Qobs & \defn & \estim{\Exs} \left[\varfun{X}{\eparams}
X_{\vnots{\svert}} X_{\vnots{\svert}}^T \right] \; = \;
\frac{1}{\numobs} \sum_{i=1}^\numobs \varfun{\myxsam{i}}{\eparams}
\myxsam{i}_{\vnots{\svert}} (\myxsam{i}_{\vnots{\svert}})^T
\end{eqnarray}
then the growth condition~\eqref{EqnGrowthCondition} and choice of
$\regpar_\numobs$ from Theorem~\ref{ThmMain} are sufficient to ensure
that the graph is recovered with high probability.  Interestingly, our
analysis shows that if the conditions are imposed directly on the
sample Fisher information matrices and $\thetamin = \Theta(1)$, then
the weaker growth condition $\numobs = \Omega(\degmax^2 \log(\pdim))$
suffices for asymptotically exact graph recovery.

The second part of the analysis, provided in
Section~\ref{SecPopulation}, is devoted to showing that under the
specified growth conditions (\assgro), imposing incoherence and
dependence assumptions on the \emph{population version} of the Fisher
information $\Qstar$ guarantees (with high probability) that analogous
conditions hold for the sample quantities $\Qobs$.  While it follows
immediately from the law of large numbers that the empirical Fisher
information $\Qobs_{AA}$ converges to the population version
$\Qstar_{AA}$ for any \emph{fixed} subset $A$, the delicacy is that we
require controlling this convergence over subsets of increasing size.
The analysis therefore requires some large-deviations bounds, so as to
provide exponential control on the rates of convergence.

\subsection{Primal-dual witness for graph recovery}
\label{SecPrimalDual}

At a high-level, at the core of our proof lies the notion of a
primal-dual witness.  In particular, we explicitly construct an
optimal \emph{primal-dual pair}, namely, a primal solution
$\estim{\eparam}$, along with an associated subgradient vector $\zhat$
(which can be interpreted as a dual solution), such that the
Karush-Kuhn-Tucker (KKT) conditions associated with the convex
program~\eqref{EqnRegLikeTwo} are satisfied.  Moreover, we show that
under the stated assumptions on $(\numobs, \pdim, \degmax)$, the
primal-dual pair $(\estim{\eparam}, \zhat)$ can be constructed such
that they act as a \emph{witness}---that is, a certificate
guaranteeing that the method correctly recovers the graph structure.

Let us write the convex program~\eqref{EqnRegLikeTwo} in the form
\begin{equation}
\label{EqnRegLikeThree}
\min \limits_{\eparam_{\vnots{\svert}} \in \real^{\pdim -1} } \left \{
\neglog(\eparam; \Data) + \regpar_\numobs
\|\eparam_{\vnots{\svert}}\|_1 \right \},
\end{equation}
where 
\begin{equation}
\neglog(\eparam; \Data) \, = \, \neglog(\eparam) \, = \,
  \frac{1}{\numobs} \sum_{i=1}^\numobs \Logfun{\eparam}{\myxsam{i}} -
  \sum_{u \in \vertex \bk \svert} \eparam_{\svert u} \muhat_{\svert u}
\end{equation}
is the negative log likelihood associated with the logistic regression
model.  The KKT conditions associated with this model can be expressed
as follows
\begin{eqnarray}
\label{EqnZeroGrad}
\nabla \neglog(\estim{\eparam}) + \regpar_\numobs \dualvec & = & 0
\end{eqnarray}
where the dual or subgradient vector $\dualvec \in \real^{\pdim -1}$
satisfies the properties
\begin{equation}
\label{EqnDualvecProp}
\dualvec_{\svert t} = \sign(\estim{\eparam}_{\svert t}) \quad \mbox{if
$\estim{\eparam}_i \neq 0$, and} \quad |\dualvec_{\svert t}| \leq 1
\quad \mbox{otherwise}.
\end{equation}
One way to understand this vector $\dualvec$ is as a subgradient,
meaning an element of the subdifferential of the $\ell_1$-norm (see
Rockafellar,~\citeyear{Rockafellar}).  An alternative interpretation is
based on the constrained equivalent to
problem~\eqref{EqnRegLikeThree}, involving the constraint
$\|\eparam\|_1 \leq C(\regpar_\numobs)$.  This $\ell_1$-constraint is
equivalent to the family of constraints $\signvec^T \eparam \leq C$,
where the vector $\signvec \in \{-1,+1\}^{\pdim -1}$ ranges over all
possible sign vectors.  In this formulation, the optimal dual vector
is simply the conic combination 
\begin{equation}
\label{EqnConvexComb}
\regpar_\numobs \dualvec = \sum_{\signvec \in \{-1,+1\}^{\pdim-1}}
\alpha^*_{\psignvec} \signvec,
\end{equation}
where $\alpha^*_\psignvec \geq 0$ is the Lagrange multiplier
associated with the constraint $\signvec^T \eparam \leq C$.

The KKT conditions~\eqref{EqnZeroGrad} and~\eqref{EqnDualvecProp} must
be satisfied by any optimal pair $(\estim{\eparam}, \dualvec)$ to the
convex program~\eqref{EqnRegLikeThree}.  In order for this primal-dual
pair to correctly specify the graph structure, we require furthermore
that the following properties are satisfied:
\begin{subequations}
\begin{align}
\label{EqnCorrectSign}
& \sign(\dualvec_{\svert t}) =  \sign(\eparam^*_{\svert t}) \qquad
\mbox{for all $(\svert, t) \in \Nset \defn \{(\svert, t) \, \mid \, t
\in \N(\svert) \}$, $\quad$ and} \\
\label{EqnNoFalsePos}
& \mythetahat_\Nsetcomp  =  0 \quad \mbox{where $\Nsetcomp \defn
\{(\svert, u) \, \mid \, (\svert, u) \notin \edge \}$}.
\end{align}
\end{subequations}

We now construct our \emph{witness} pair $(\mythetahat,\dualvec)$ as
follows.  First, we set $\mythetahat_{\Nset}$ as the minimizer of the
partial penalized likelihood, 
\begin{eqnarray}
\label{thetaSCon}
\mythetahat_{\Nset} = \arg\min_{(\eparam_{\Nset}, 0) \in \real^{\pdim
    -1}} \left\{\neglog(\eparam; \Data) +
\regpar_\numobs\|\eparam_{\Nset}\|_1 \right\},
\end{eqnarray}
and set $\dualvec_{\Nset} = \sign(\mythetahat_{\Nset})$.  We then set
$\mythetahat_\Nsetcomp = 0$ so that condition~\eqref{EqnNoFalsePos}
holds. Finally, we obtain $\dualvec_{\Nsetcomp}$ from
equation~\eqref{EqnZeroGrad} by plugging in the values of
$\mythetahat$ and $\dualvec_{\Nset}$. Thus, our construction satisfies
conditions~\eqref{EqnNoFalsePos} and \eqref{EqnZeroGrad}. The
remainder of the analysis consists of showing that our
conditions on $(\numobs, \pdim, \degmax)$ imply that, with high-probability,
the remaining conditions~\eqref{EqnCorrectSign}
and~\eqref{EqnDualvecProp} are satisfied.

This strategy is justified by the following lemma, which provides sufficient conditions for
shared sparsity and uniqueness of the optimal solution:
\blems[Shared sparsity and uniqueness]
If there exists an optimal primal solution $\estim{\eparam}$ with
associated optimal dual vector $\dualvec$ such that
$\|\dualvec_{\Nsetcomp}\|_\infty < 1$, then any optimal primal
solution $\wtil{\eparam}$ must have $\wtil{\eparam}_{\Nsetcomp} = 0$.
Moreover, if the Hessian sub-matrix $[\nabla^2
\neglog(\estim{\eparam})]_{\Sset \Sset}\succ 0$, then
$\estim{\eparam}$ is the unique optimal solution.
\elems
\spro

By Lagrangian duality, the penalized problem~\eqref{EqnRegLikeThree}
can be written as an equivalent constrained optimization problem over
the ball $\|\eparam\|_1 \leq C(\regpar_\numobs)$, for some constant
$C(\regpar_\numobs) < +\infty$.  Since the Lagrange multiplier
associated with this constraint---namely, $\regpar_\numobs$---is
strictly positive, the constraint is active at any optimal solution,
so that $\|\eparam\|_1$ is constant across all optimal solutions.
Consider the representation of $\dualvec$ as the convex
combination~\eqref{EqnConvexComb} of sign vectors $\signvec \in
\{-1,+1\}^{\pdim -1}$, where the weights $\alpha^*_\psignvec$ are
non-negative and sum to one.  Since $\alpha^*$ is an optimal vector of
Lagrange multipliers for the optimal primal solution
$\estim{\eparam}$, it follows~\citep{Bertsekas_nonlin} that any other
optimal primal solution $\wtil{\eparam}$ must minimize the associated
Lagrangian (i.e., satisfy equation~\eqref{EqnZeroGrad}), and moreover
must satisfy the complementary slackness conditions
$\alpha^*_\psignvec [\signvec^T \eparam - C] = 0$ for all sign vectors
$v$.  But these conditions imply that $\dualvec^T \wtil{\eparam} = C =
\|\wtil{\eparam}\|_1$, which cannot occur if $\wtil{\eparam}_j \neq 0$
for some index $j$ for which $|\dualvec_j| < 1$.  We thus conclude
that $\wtil{\eparam}_{\Nsetcomp} = 0$ for all optimal primal
solutions.

Finally, given that all optimal solutions satisfy $\eparam_{\Nsetcomp}
= 0$, we may consider the restricted optimization problem subject to
this set of constraints.  If the principal submatrix of the Hessian is
positive definite, then this sub-problem is strictly convex, so that
the optimal solution must be unique.
\fpro

In our primal-dual witness proof, we exploit this lemma by
constructing a primal-dual pair $(\estim{\eparam}, \dualvec)$ such
that $\|\dualvec_{S^c}\|_\infty < 1$.  Moreover, under the conditions of
Theorem~\ref{ThmMain}, we prove that the sub-matrix of the sample
Fisher information matrix is strictly positive definite with high
probability, so that the primal solution $\estim{\eparam}$ is
guaranteed to be unique.

\section{Analysis under sample Fisher matrix assumptions}
\label{SecFixedDesign}

We begin by establishing model selection consistency when assumptions
are imposed directly on the sample Fisher matrix $\Qobs$, as opposed
to on the population matrix $\Qstar$, as in Theorem~\ref{ThmMain}.
In particular, we define the ``good event'' 
\begin{eqnarray}
\GoodEvent(\Data) & \defn & \left\{ \Data \, \mid \, \mbox{$\Qobs =
\estim{\Exs}[\nabla^2 \neglog(\eparams)]$ satisfies $\assdep$ and
$\assinc$} \right \}.
\end{eqnarray}
We then state the following
\bprops[Consistency for fixed design]
\label{PropFixed}
If $\numobs > L \degmax^2 \log(\pdim)$ for a suitably large
constant $L$, and the minimum value $\thetamin$ decays no faster than
$\order(1/\sqrt{\degmax})$, then for the regularization sequence
$\regpar_\numobs = 2\sqrt{\frac{\log \pdim}{\numobs}}$, the estimated
graph $\estim{\graph}(\regpar_\numobs)$ obtained by neighborhood
logistic regression satisfies
\begin{align}
\label{EqnSuccessTwo}
\mprob[\estim{\graph}(\regpar_\numobs) = \graph(\pdim, \degmax) \,
\mid \, \GoodEvent(\Data)] & = & \order \left(\exp\left(-\numobs
\regpar_\numobs^2 - \log(\pdim) \right)\right) \; \rightarrow \; 0.
\end{align}
\eprops
Loosely stated, this result guarantees that if the sample Fisher
information matrix is ``good'', then the conditional probability of
successful graph recovery converges to zero at the specified rate.
The remainder of this section is devoted to the proof of
Proposition~\ref{PropFixed}.

\subsection{Key technical results}
\label{SecTechnical}

We begin with statements of some key technical lemmas that are central to
our main argument, with their proofs deferred to
Appendix~\ref{AppTechnical}. The central object is the following
expansion, obtained by re-writing the zero-subgradient condition as
\begin{eqnarray}
\label{EqnTemp}
\nabla \neglog(\eparamhat; \Data) - \nabla \neglog (\eparams; \Data) &
= & \obsnoise - \regpar_\numobs \dualvec,
\end{eqnarray}
where we have introduced the short-hand notation $\obsnoise$ for the
$(\pdim-1)$-vector
\begin{eqnarray*}
\lefteqn{\obsnoise \defn -\nabla \neglog(\eparams; \Data) \; = \;} && \\
&& -\frac{1}{\numobs} \sum_{i=1}^\numobs \myxsam{i}_{\bk \svert}
\left\{\myxsam{i}_\svert - \frac{\exp(\sum_{t \in \vertex \bk \svert}
  \eparam^*_{\svert t} \myxsam{i}_t) -\exp(-[\sum_{t \in \vertex \bk
      \svert} \eparam^*_{\svert t} \myxsam{i}_t]) } {\exp(\sum_{t \in
    \vertex \bk \svert} \eparam^*_{\svert t} \myxsam{i}_t)
  +\exp(-[\sum_{t \in \vertex \bk \svert} \eparam^*_{\svert t}
    \myxsam{i}_t]) } \right\} \\
& = & -\frac{1}{\numobs} \sum_{i=1}^\numobs \myxsam{i}_{\bk \svert}
\left\{\myxsam{i}_\svert - \mprob_{\eparams}[x_\svert = 1 \, \mid \,
\myxsam{i}_{\bk \svert}] + \mprob_{\eparams}[x_\svert = -1 \, \mid \,
\myxsam{i}_{\bk \svert}] \right\}.
\end{eqnarray*}
For future reference, note that $\Exs_{\eparams}[\obsnoise] = 0$.
Next, applying the mean-value theorem coordinate-wise to the expansion~\eqref{EqnTemp}
yields
\begin{eqnarray}
\label{EqnTay}
\nabla^2 \neglog(\eparams; \Data) \; [\estim{\eparam} - \eparams] & =
& \obsnoise -\regpar_\numobs \dualvec + \Rem^\numobs,
\end{eqnarray}
where the remainder term takes the form
\begin{eqnarray}
\label{EqnRemainder}
\Rem^\numobs_j & = & \left[\nabla^2 \neglog(\eparambar^{(j)}; \Data) -
\nabla^2 \neglog(\eparams; \Data) \right]_j^T \, (\eparamhat - \eparams),
\end{eqnarray}
with $\eparambar^{(j)}$ a parameter vector on the line between $\eparams$
and $\eparamhat$, and with $[\cdot]_j^T$ denoting the $j$th row of the matrix.  
The following lemma addresses the behavior of the
term $\obsnoise$ in this expansion:
\blems
\label{LemBernOne}
If $\numobs \regpar_\numobs^2> \log(\mdim)$, then for the specified
mutual incoherence parameter $\mutinc \in (0,1]$, we have
\begin{eqnarray}
\Prob\left(\frac{2-\mutinc}{\regpar_\numobs} \|\obsnoise\|_\infty \geq
\frac{\mutinc}{4}\right) & = & \order\left(\exp\left(- \numobs
\regpar_\numobs^2 + \log(\mdim) \right)\right) \rightarrow 0.
\end{eqnarray}
\elems
\noindent See Appendix~\ref{AppBernOne} for the proof of this claim.

The following lemma establishes that the sub-vector
$\widehat{\eparam}_\Nset$ is an $\ell_2$-consistent estimate of the
true sub-vector $\eparam^*_\Nset$:
\blems[$\ell_2$-consistency of primal subvector]
\label{Leml2cons}
If $\regpar_\numobs \degmax \leq \frac{\Cmin^2}{10\Dmax}$, then as
$\numobs \rightarrow +\infty$, we have
\begin{equation}
\label{Eqnl2cons}
\|\hat{\eparam}_{\Nset} - \tparam_{\Nset}\|_{2} =
\order_p\left(\sqrt{\degmax}\lambda_{\numobs}\right) \, \rightarrow \, 0.
\end{equation}
\elems
\noindent See Appendix~\ref{Appl2cons} for the proof of this claim.

Our final technical lemma provides control on the the remainder
term~\eqref{EqnRemainder}:
\blems
\label{LemTayRem}
If $\numobs \regpar^2_\numobs > \log(\mdim)$ and $\degmax
\regpar_\numobs$ is sufficiently small, then for mutual incoherence
parameter $\mutinc \in (0,1]$, we have
\begin{eqnarray}
\Prob\left(\frac{2 - \mutinc}{\regpar_\numobs} \|\Rem^\numobs\|_\infty
\geq \frac{\mutinc}{4}\right) & = & \order(\exp\left(-\numobs
\regpar_\numobs^2 + \log(\mdim) \right)) \; \rightarrow \; 0.
\end{eqnarray}
\elems
\noindent See Appendix~\ref{AppTayRem} for the proof of this claim.

%%%%%%%%%%%%%%%%%%%%%%%%%%%%%%%%%%%%%%%%%%%%%%%%%%%%%%%%%%%%%%%%%%%%%%%%%%%%%%%%

\subsection{Proof of Proposition~\ref{PropFixed}}
Using these lemmas, we can now complete the proof of
Proposition~\ref{PropFixed}.  Recalling our shorthand \mbox{$\Qobs =
\nabla^2_\theta \neglog(\eparams; \Data)$,} we re-write
condition~\eqref{EqnTay} in block form as:
\begin{subequations}
\label{EqnBlock}
\begin{eqnarray}
\label{EqnBlockA}
\Qobs_{\Sbar \Sset} \; [\thetahat_\Nset -\thetastar_\Nset] & = &
\obsnoise_{\Sbar} -\regpar_\numobs \dualvec_\Sbar +
\Rem^\numobs_{\Sbar}, \\
\label{EqnBlockB}
\Qobs_{\Sset \Sset} \; [\thetahat_\Nset -\thetastar_\Nset
] & = & \obsnoise_\Sset -\regpar_\numobs
\dualvec_\Sset + \Rem^\numobs_\Sset.
\end{eqnarray}
\end{subequations}
Since the matrix $\Qobs_{\Sset \Sset}$ is invertible by assumption,
the conditions~\eqref{EqnBlock} can be re-written as
\begin{eqnarray}
\Qobs_{\Sbar \Sset} \; (\Qobs_{\Sset \Sset})^{-1} \;
\left[\obsnoise_\Sset -\regpar_\numobs \dualvec_\Sset +
\Rem^\numobs_\Sset \right] & = & \obsnoise_{\Sbar} -\regpar_\numobs
\dualvec_\Sbar + \Rem^\numobs_{\Sbar}.
\end{eqnarray}
Rearranging yields the condition
\begin{align}
\label{EqnDefnDualsub}
\left[\obsnoise_{\Sbar} - \Rem^\numobs_{\Sbar} \right] -
\Qobs_{\Sbar \Sset} \, (\Qobs_{\Sset \Sset})^{-1}
\left[\obsnoise_\Sset - \Rem^\numobs_\Sset \right] +
\regpar_\numobs \Qobs_{\Sbar \Sset} \, (\Qobs_{\Sset \Sset})^{-1}
\dualvec_\Sset & =  \regpar_\numobs \dualvec_\Sbar.
\end{align}

\myparagraph{Strict dual feasibility} We now demonstrate that for the
dual sub-vector $\dualvec_\Sbar$ defined by
equation~\eqref{EqnDefnDualsub}, we have $\|\dualvec_\Sbar\|_\infty <
1$.  Using the triangle inequality and the mutual incoherence
bound~\eqref{EqnDefnMutualInco}, we have that
\begin{eqnarray}
\label{EqnTriangle}
\|\dualvec_\Sbar\|_\infty & \leq & \matsnorm{\Qobs_{\Sbar \Sset} \,
(\Qobs_{\Sset \Sset})^{-1}}{\infty}
\left[\frac{\|\obsnoise_\Sset\|_\infty}{\regpar_\numobs} +
\frac{\|\Rem^\numobs_\Sset\|_\infty}{\regpar_\numobs} + 1\right]\\
&&\qquad +
\frac{\|\Rem^\numobs_\Sbar\|_\infty}{\regpar_\numobs} +
\frac{\|\obsnoise_\Sbar\|_\infty}{\regpar_\numobs} \nonumber \\
& \leq & (1-\mutinc) + (2-\mutinc)
\left[\frac{\|\Rem^\numobs\|_\infty}{\regpar_\numobs} +
\frac{\|\obsnoise\|_\infty}{\regpar_\numobs} \right]
\end{eqnarray}
Next, applying Lemmas~\ref{LemBernOne} and~\ref{LemTayRem}, we have
\begin{eqnarray*}
\|\dualvec_\Sbar\|_\infty & \leq & (1-\mutinc) + \frac{\mutinc}{4} +
\frac{\mutinc}{4} \; = \; 1 - \frac{\mutinc}{2}.
\end{eqnarray*}
with probability converging to one.

\myparagraph{Correct sign recovery:} We next show that our primal
sub-vector $\estim{\eparam}_\Sset$ defined by
equation~\eqref{thetaSCon} satisfies sign consistency, meaning,
$\sgn(\estim{\eparam}_{\Sset}) = \sgn(\eparam^*_{\Sset})$. In order to establish this, 
it suffices to show that
\begin{eqnarray*}
\left \|\eparam_{\Nset} - \tparam_{\Nset} \right \|_\infty &
\leq & \frac{\thetamin}{2}
\end{eqnarray*}
where we recall the notation $\thetamin \defn \min_{(\svert, t) \in
\edge} |\thetastar_{\svert t}|$.  From Lemma~\ref{Leml2cons}, we have
$\|\eparam_{\Nset} - \tparam_{\Nset}\|_{2} =
\order_{p}(\sqrt{d}\regpar_\numobs)$, so that \beqn
\frac{2}{\thetamin}\|\eparam_{\Nset} - \tparam_{\Nset}\|_\infty &\le&
\frac{2}{\thetamin}\|\eparam_{\Nset} - \tparam_{\Nset}\|_2\\ &=&
\order\left(\frac{\sqrt{d} \lambda_{n}}{\thetamin}\right) \eeqn

Since $\thetamin$ decays no faster than $\Theta(1/\sqrt{\degmax})$,
the right-hand side is upper bounded by $\order(\regpar_\numobs
\degmax)$, which can be made smaller than $1$ by choosing
$\regpar_\numobs$ sufficiently small, as asserted in
Proposition~\ref{PropFixed}.

\section{Uniform convergence of sample information matrices}
\label{SecPopulation}

In this section, we complete the proof of Theorem~\ref{ThmMain} by
showing that if the dependency ($\assdep$) and incoherence
($\assinc$) assumptions are imposed on the \emph{population} Fisher
information matrix then under the specified scaling of $(\numobs,
\pdim, \degmax)$, analogous bounds hold for the \emph{sample} Fisher
information matrices with probability converging to one.  These
results are not immediate consequences of classical random matrix
theory (e.g.,~\cite{DavSza01}), since the elements of $\Qobs$ are
highly dependent.

Recall the definitions
\begin{equation}
\label{EqnRecall}
\Qstar \defn \Exs_{\eparams} \left[\varfun{X}{\eparams}
X_{\vnots{\svert}} X_{\vnots{\svert}}^T \right], \quad \mbox{and}
\quad \Qobs \defn \estim{\Exs}\left[\varfun{X}{\eparams}
X_{\vnots{\svert}} X_{\vnots{\svert}}^T \right],
\end{equation}
where $\Exs_\eparams$ denotes the population expectation, and
$\estim{\Exs}$ denotes the empirical expectation, and the variance
function $\eta$ was defined previously equation~\eqref{EqnDefnVarfun}.  The
following lemma asserts the eigenvalue bounds in Assumption $\assdep$
hold with high probability for sample covariance matrices:
\blems
\label{LemConcentrateEig}
Suppose that assumption $\assdep$ holds for the population matrix
$\Qstar$ and $\Exs_{\eparams}[X X^T]$.  For any $\delta > 0$ and some
fixed constants $A$ and $B$, we have
\begin{subequations}
\begin{align}
\label{EqnEigConcMax}
\Prob\left[\myeigmax\big[\frac{1}{\numobs} \sum_{i=1}^\numobs
\myxsam{i}_{\vnots{\svert}} (\myxsam{i}_{\vnots{\svert}})^T\big] \geq \Dmax
+ \delta\right] & \leq  2 \exp \left(-A\frac{\delta^2 \numobs}{ \degmax^2}
+ B \log(\degmax)\right). \\
\label{EqnEigConcMin}
\Prob[\myeigmin(\Qobs_{\Sset\Sset}) \leq \Cmin - \delta] & \leq  2
\exp \left(-A\frac{\delta^2 \numobs}{ \degmax^2} + B \log(\degmax)
\right).
\end{align}
\end{subequations}
\elems

The following result is the analog for the incoherence assumption
($\assinc$), showing that the scaling of $(\numobs, \mdim, \degmax)$
given in Theorem~\ref{ThmMain} guarantees that population incoherence
implies sample incoherence.
\blems 
\label{LemConcentrateInco}
If the population covariance satisfies a mutual incoherence
condition~\eqref{EqnDefnMutualInco} with parameter $\mutinc \in (0,1]$
as in Assumption $\assinc$, then the sample matrix satisfies an
analogous version, with high probability, in the sense that
\begin{eqnarray}
\label{EqnSampleMutinc}
\mprob \left[\matnorm{\Qobs_{\Sbar \Sset} (\Qobs_{\Sset
\Sset})^{-1}}{\infty}{\infty} \geq 1- \frac{\mutinc}{2} \right] & \leq
& \exp \left(-\unicon \frac{\numobs}{\degmax^3} + \log(\pdim) \right).
\end{eqnarray}
\elems

Proofs of these two lemmas are provided in the following sections.
Before proceeding, we begin by taking note of a simple bound to be
used repeatedly throughout our arguments.  By definition of the
matrices $\Qobs(\eparam)$ and $\Qmat(\eparam)$ (see
equations~\eqref{EqnDefnSampleFisher} and~\eqref{EqnDefnQstarExp}),
the $(j,k)^{th}$ element of the difference matrix $\Qobs(\eparam) -
\Qmat(\eparam)$ can be written as an i.i.d. sum of the form $Z_{jk} =
\frac{1}{\numobs}\sum_{i=1}^\numobs \Zi_{jk}$, where each $\Zi_{jk}$
is zero-mean and bounded (in particular, $|\Zi_{jk}| \leq 4$).  By the
Azuma-Hoeffding bound~\citep{Hoeffding63}, for any indices $j,k = 1,
\ldots, \degmax$ and for any $\epsilon > 0$, we have
\begin{eqnarray}
\label{EqnAzuma}
\mprob[(\myZ_{jk})^2 \geq \epsilon^2 ] \; = \;
\mprob[|\frac{1}{\numobs} \sum_{i=1}^\numobs \Zi_{jk}| \geq \epsilon]
& \leq & 2 \exp \left( -\frac{\epsilon^2 \numobs}{32} \right).
\end{eqnarray}
So as to simplify notation, throughout this section, we use $\unicon$
to denote a universal positive constant, independent of $(\numobs,
\pdim, \degmax)$.  Note that the precise value and meaning of
$\unicon$ may differ from line to line.

\subsection{Proof of Lemma~\ref{LemConcentrateEig}}

By the Courant-Fischer variational representation~\citep{Horn85}, we
have
\begin{eqnarray*}
\myeigmin(\Qmat_{\Sset \Sset}) & = & \min_{\|x\|_2 = 1} x^T
\Qmat_{\Sset \Sset} x \\
& = & \min_{\|x\|_2 = 1} \left \{ x^T \Qobs_{\Sset \Sset} x +
x^T(\Qmat_{\Sset \Sset} - \Qobs_{\Sset \Sset})x \right \} \\
& \leq & y^T \Qobs_{\Sset \Sset} y + y^T(\Qmat_{\Sset \Sset} -
\Qobs_{\Sset \Sset}) y,
\end{eqnarray*}
where $y \in \real^{\degmax}$ is a unit-norm minimal eigenvector of
$\Qobs_{\Sset \Sset}$.  Therefore, we have
\begin{eqnarray*}
\myeigmin(\Qobs_{\Sset \Sset}) & \geq & \myeigmin(\Qmat_{\Sset \Sset})
- \matsnorm{\Qmat_{\Sset \Sset}- \Qobs_{\Sset \Sset}}{2} \; \geq \;
\Cmin - \matsnorm{\Qmat_{\Sset \Sset}- \Qobs_{\Sset \Sset}}{2}.
\end{eqnarray*}
Hence it suffices to obtain a bound on the spectral norm
$\matsnorm{\Qmat_{\Sset \Sset}- \Qobs_{\Sset \Sset}}{2}$.  Observe
that
\begin{eqnarray*}
\matnorm{\Qobs_{\Sset \Sset} - \Qmat_{\Sset \Sset}}{2}{2} & \leq &
    \big[ \sum_{j=1}^\degmax \sum_{k=1}^\degmax (\myZ_{jk})^2
    \big]^{\frac{1}{2}}.
\end{eqnarray*}
Setting $\epsilon^2 = \delta^2/\degmax^2$ in equation~\eqref{EqnAzuma}
and applying the union bound over the $\degmax^2$ index pairs $(j,k)$ then
yields
\begin{eqnarray}
\label{EqnUseful}
\mprob[\matsnorm{\Qobs_{\Sset \Sset} -\Qmat_{\Sset \Sset}}{2} \geq
  \delta ] & \leq & 2 \exp \left(-\unicon \frac{\delta^2
  \numobs}{\degmax^2} + 2 \log(\degmax)\right).
\end{eqnarray}

Similarly, we have
\begin{align*}
& \Prob[\myeigmax(\frac{1}{\numobs} \sum_{i=1}^\numobs
\myxsam{i}_{\vnots{\svert}} (\myxsam{i}_{\vnots{\svert}})^T) \geq \Dmax]
\leq\\
&\qquad  \Prob \left[\matsnorm{ \big(\frac{1}{\numobs}
\sum_{i=1}^\numobs \myxsam{i}_{\vnots{\svert}}
(\myxsam{i}_{\vnots{\svert}})^T \big) -
\Exs_\eparams[X_{\vnots{\svert}} X_{\vnots{\svert}}^T]}{2} \geq \delta
\right],
\end{align*}
which obeys the same upper bound~\eqref{EqnUseful}, by following the
analogous argument.

%%%%%%%%%%%%%%%%%%%%%%%%%%%%%%%%%%%%%%%%%%%%%%%%%%%%%%%%%%%%%%%%%%%%%%%%%%%

\subsection{Proof of Lemma~\ref{LemConcentrateInco}}

We begin by decomposing the sample matrix as the sum $ \Qobs_{\Sbar
\Sset} (\Qobs_{\Sset \Sset})^{-1} = \Term_1 + \Term_2 + \Term_3 +
\Term_4$, where we define
\begin{subequations}
\begin{eqnarray}
\Term_1 & \defn & \Qstar_{\Sbar \Sset} \left[ (\Qobs_{\Sset
\Sset})^{-1} - (\Qstar_{\Sset \Sset})^{-1} \right]  \\
\Term_2 & \defn & \left[\Qobs_{\Sbar \Sset} - \Qstar_{\Sbar \Sset}
\right] \, (\Qstar_{\Sset \Sset})^{-1}\\ 
\Term_3 & \defn & \left[\Qobs_{\Sbar \Sset}
- \Qstar_{\Sbar \Sset}\right] \left[ (\Qobs_{\Sset \Sset})^{-1} -
(\Qstar_{\Sset \Sset})^{-1}\right] \\ 
\Term_4 & \defn & \Qstar_{\Sbar \Sset} (\Qstar_{\Sset \Sset})^{-1}
\end{eqnarray}
\end{subequations}
The fourth term is easily controlled; indeed, we have
$$\matnorm{\Term_4}{\infty}{\infty} = \matnorm{\Qstar_{\Sbar \Sset}
(\Qstar_{\Sset \Sset})^{-1}}{\infty}{\infty} \leq 1 - \mutinc$$ by the
incoherence assumption $\assinc$.  If we can show that
$\matnorm{\Term_i}{\infty}{\infty} \leq \frac{\mutinc}{6}$ for the
remaining indices $i=1,2,3$, then by our four term decomposition and
the triangle inequality, the sample version satisfies the
bound~\eqref{EqnSampleMutinc}, as claimed.  We deal with these
remaining terms using the following lemmas:
\blems\label{LemTechBounds}
For any $\delta > 0$ and constants $\unicon, \unicon'$, the following
bounds hold:
\begin{subequations}
\label{EqnTechBound}
\begin{align}
\label{EqnTechBoundA}
\mprob[\matnorm{\Qobs_{\Sbar \Sset} - \Qstar_{\Sbar
\Sset}}{\infty}{\infty} \geq \delta] & \leq  2 \exp \left( -\unicon
\frac{\numobs \, \delta^2}{\degmax^2} + \log(\degmax) + \log(\pdim -
\degmax) \right) \\
\label{EqnTechBoundB}
\mprob[\matnorm{\Qobs_{\Sset \Sset} - \Qstar_{\Sset
\Sset}}{\infty}{\infty} \geq \delta] & \leq  2 \exp \left( -\unicon
\frac{\numobs \, \delta^2}{\degmax^2} + 2 \log(\degmax) \right) \\
\label{EqnTechBoundC}
\mprob[\matnorm{(\Qobs_{\Sset \Sset})^{-1} - (\Qstar_{\Sset
\Sset})^{-1}}{\infty}{\infty} \geq \delta] & \leq  4 \exp \left(
-\unicon \frac{\numobs \, \delta^2}{\degmax^3} + \unicon' \log(\degmax)
\right).
\end{align}
\end{subequations}
\elems
\noindent See Appendix~\ref{AppTechBounds} for the proof of these
claims.

\myparagraph{Control of first term:} Turning to the first term, we
first re-factorize it as
\begin{eqnarray*}
\Term_1 & = & \Qstar_{\Sbar \Sset} (\Qstar_{\Sset \Sset})^{-1}
\left[\Qobs_{\Sset \Sset} - \Qstar_{\Sset \Sset} \right]
(\Qobs_{\Sset \Sset})^{-1},
\end{eqnarray*}
and then bound it (using the sub-multiplicative property\\
$\matsnorm{A B}{\infty} \leq \matsnorm{A}{\infty} \matsnorm{B}{\infty}$) as follows
\begin{eqnarray*}
\matnorm{\Term_1}{\infty}{\infty} & \leq & \matnorm{\Qstar_{\Sbar
\Sset} (\Qstar_{\Sset \Sset})^{-1}}{\infty}{\infty}
\matnorm{\Qobs_{\Sset \Sset} - \Qstar_{\Sset \Sset}}{\infty}{\infty}
\matnorm{(\Qobs_{\Sset \Sset})^{-1}}{\infty}{\infty} \\
& \leq & (1-\mutinc) \, \matnorm{\Qobs_{\Sset \Sset} - \Qstar_{\Sset
\Sset}}{\infty}{\infty} \left \{ \sqrt{\degmax} \; \matnorm{(\Qobs_{\Sset
\Sset})^{-1}}{2}{2} \right \},
\end{eqnarray*}
where we have used the incoherence assumption $\assinc$.  Using the
bound~\eqref{EqnEigConcMin} from Lemma~\ref{LemConcentrateEig} with $\delta =
\Cmin/2$, we have $\matnorm{(\Qobs_{\Sset \Sset})^{-1}}{2}{2} \; = \;
[\myeigmin(\Qobs_{\Sset \Sset})]^{-1} \; \leq \; \frac{2}{\Cmin}$ with
probability greater than $1-\exp \left(-\unicon \numobs/ \degmax^2+ 2
\log(\degmax)\right)$.  Next, applying the bound~\eqref{EqnTechBoundB}
with $\delta = c/\sqrt{\degmax}$, we conclude that with probability
greater than $1-2 \exp \left( -\unicon \numobs c^2/\degmax^3 +
\log(\degmax) \right)$, we have
\begin{equation*}
\matnorm{\Qobs_{\Sset \Sset} - \Qstar_{\Sset
\Sset}}{\infty}{\infty} \; \leq \;  c/\sqrt{\degmax}.
\end{equation*}
By choosing the constant $c > 0$ sufficiently small, we are guaranteed
that
\begin{eqnarray}
\mprob[ \matnorm{\Term_1}{\infty}{\infty} \geq \mutinc/6] & \leq & 2
 \exp \left( -\unicon \frac{\numobs c^2}{\degmax^3} + \log(\degmax)
 \right).
\end{eqnarray}

\myparagraph{Control of second term:} To bound $\Term_2$, we first
write
\begin{eqnarray*}
\matnorm{\Term_2}{\infty}{\infty} & \leq & \sqrt{\degmax}
\matnorm{(\Qstar_{\Sset \Sset})^{-1}}{2}{2} \; \matnorm{\Qobs_{\Sbar
\Sset} - \Qstar_{\Sbar \Sset}}{\infty}{\infty} \\
& \leq & \frac{\sqrt{\degmax}}{\Cmin} \; \matnorm{\Qobs_{\Sbar \Sset}
  - \Qstar_{\Sbar \Sset}}{\infty}{\infty}
\end{eqnarray*}
We then apply bound~\eqref{EqnTechBoundA} with $\delta =
\frac{\mutinc}{3} \, \frac{\Cmin}{\sqrt{\degmax}}$ to conclude that
\begin{eqnarray}
\mprob[\matnorm{\Term_2}{\infty}{\infty} \geq \mutinc/3] & \leq & 2
\exp \left( -\unicon \frac{\numobs}{\degmax^3} + \log(\pdim - \degmax)
\right).
\end{eqnarray}

\myparagraph{Control of third term:} Finally, in order to bound the
third term $\Term_3$, we apply the bounds~\eqref{EqnTechBoundA}
and~\eqref{EqnTechBoundB}, both with $\delta = \sqrt{\mutinc/3}$, and
use the fact that $\log(\degmax) \leq \log(\pdim - \degmax)$ to
conclude that
\begin{eqnarray}
\mprob[ \matnorm{\Term_3}{\infty}{} \geq \mutinc/3] & \leq & 4 \exp
 \left( -\unicon \frac{\numobs}{\degmax^3} + \log(\pdim - \degmax)
 \right).
\end{eqnarray}

Putting together all of the pieces, we conclude that 
\begin{eqnarray*}
\mprob[ \matnorm{\Qobs_{\Sbar \Sset} (\Qobs_{\Sset
\Sset})^{-1}}{\infty}{} \geq 1 - \mutinc/2] & = & \order
\left(\exp(-\unicon \frac{\numobs}{\degmax^3} + \log(\pdim) ) \right).
\end{eqnarray*}
as claimed.

\section{Experimental results}
\label{SecExperiments}

We now describe experimental results that illustrate some
consequences of Theorem~\ref{ThmMain}, for various types of graphs and
scalings of $(\numobs, \pdim, \degmax)$.  In all cases, we solved the
$\ell_1$-regularized logistic regression using special purpose
interior-point code developed and described by~\cite{KohKimBoy07}.

\begin{figure}
\begin{center}
\begin{tabular}{ccccc}
\widgraph{0.28\textwidth}{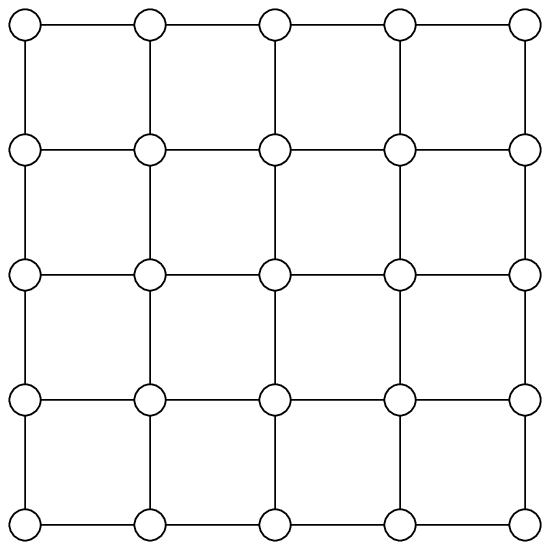} & &
\raisebox{.1in}{\widgraph{0.25\textwidth}{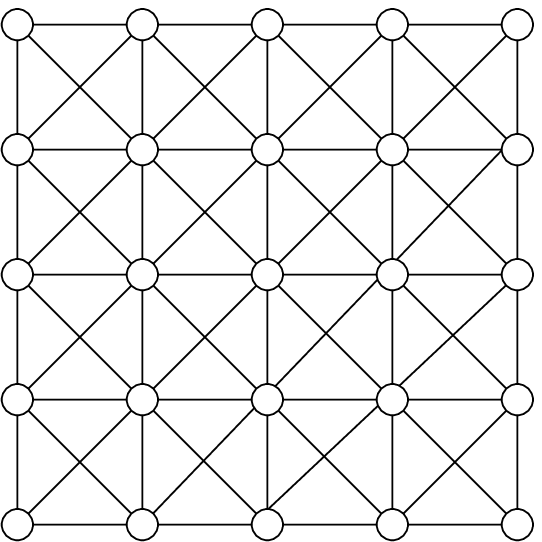}} & &
\raisebox{.1in}{\widgraph{0.25\textwidth}{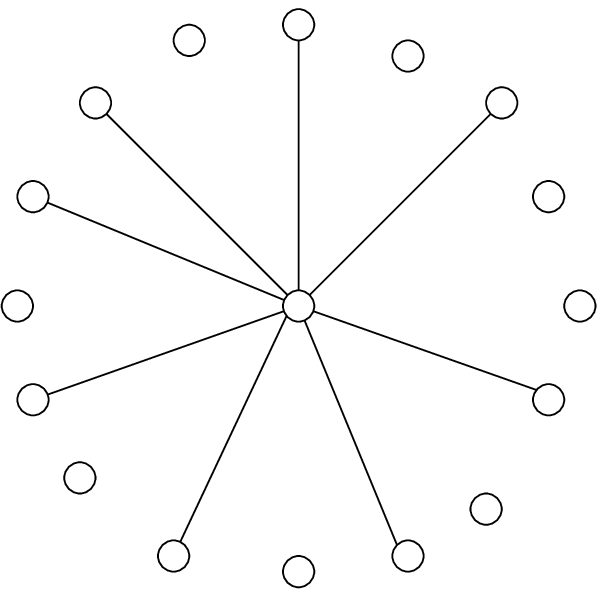}} \\
(a) & & (b) & & (c)
\end{tabular}
\caption{Illustrations of different graph classes used in simulations.
(a) Four-nearest neighbor grid ($\degmax = 4$).  (b) Eight-nearest
neighbor grid ($\degmax = 8$).  (c) Star-shaped graph ($\degmax =
\Theta(\pdim)$, or $\degmax = \Theta(\log(\pdim))$).}
\label{FigGraphs}
\end{center}
\end{figure}

We performed experiments for three different classes of graphs:
four-nearest neighbor lattices, (b) eight-nearest neighbor lattices,
and (c) star-shaped graphs, as illustrated in Figure~\ref{FigGraphs}.
Given a distribution $\mprob_{\eparams}$ of the Ising
form~\eqref{EqnIsing}, we generated random data sets $\{\myxsam{1},
\ldots, \myxsam{\numobs} \}$ by Gibbs sampling for the lattice models,
and by exact sampling for the star graph.  For a given graph class and
edge strength $\omega > 0$, we examined the performance of models with
\emph{mixed couplings}, meaning $\eparam^*_{st} = \pm \omega$ with
equal probability, or with \emph{positive couplings}, meaning that
$\eparam^*_{st} = \omega$ for all edges $(s,t)$.  In all cases, we set
the regularization parameter as $\regpar_\numobs =
\Theta(\sqrt{\frac{\log \pdim}{\numobs}})$.  Above the threshold
sample size $\numobs$ predicted by Theorem~\ref{ThmMain}, this choice
ensured correct model selection with high probability, consistent with
the theoretical prediction.  For any given graph and coupling type, we
performed simulations for sample sizes $\numobs$ scaling as $\numobs =
10 \contpar \degmax \log(\pdim)$, where the control parameter
$\contpar$ ranged from $0.1$ to upwards of $2$, depending on the graph
type.

\begin{figure}
\begin{center}
\begin{tabular}{cc}
\widgraph{0.45\textwidth}{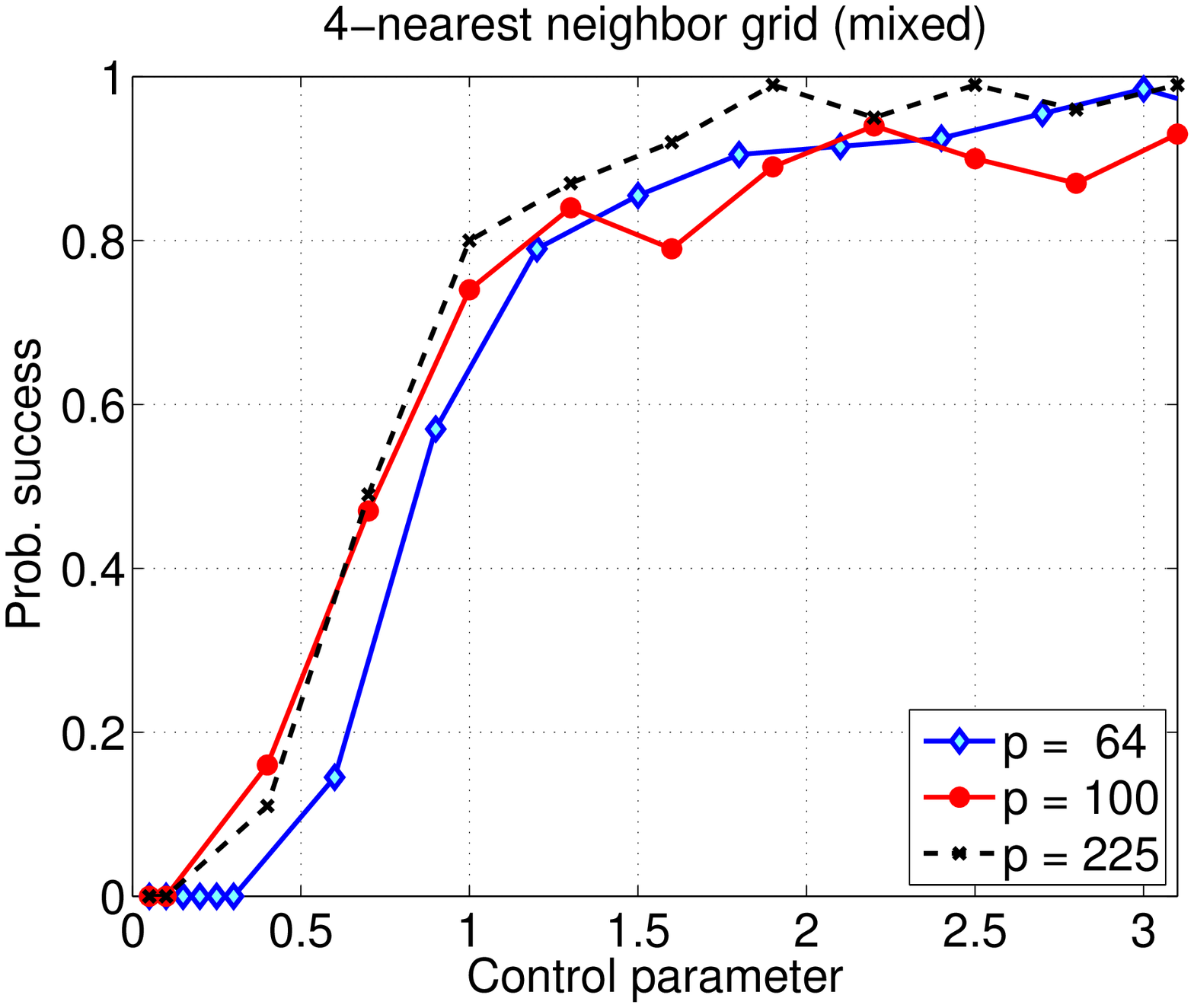} &
\widgraph{0.45\textwidth}{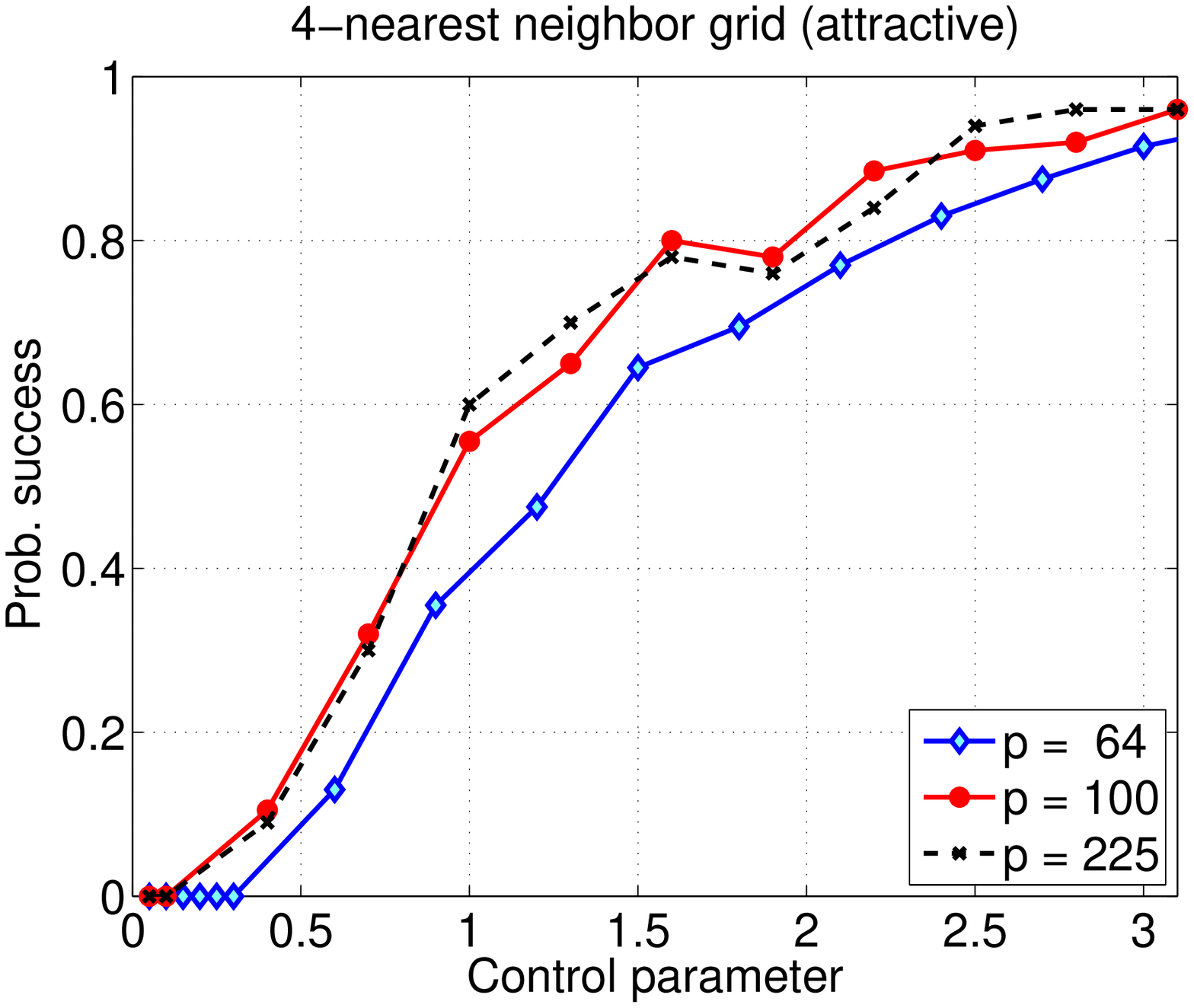} \\
(a) & (b)
\end{tabular}
\caption{Plots of success probability $\mprob[\hatNsign(\svert) =
\N(\svert), \forall r]$ versus the control parameter \mbox{$\contpar(\numobs,
\pdim, \degmax) = \numobs/[10 \degmax \log(\pdim)]$} for Ising models
on 2-D grids with four nearest-neighbor interactions ($\degmax = 4$).
(a) Randomly chosen mixed sign couplings $\eparam^*_{st} = \pm 0.50$.
(b) All positive couplings $\eparam^*_{st} = 0.50$. }
\label{FigIsingFour}
\end{center}
\end{figure}

\begin{figure}
\begin{center}
\begin{tabular}{cc}
\widgraph{0.45\textwidth}{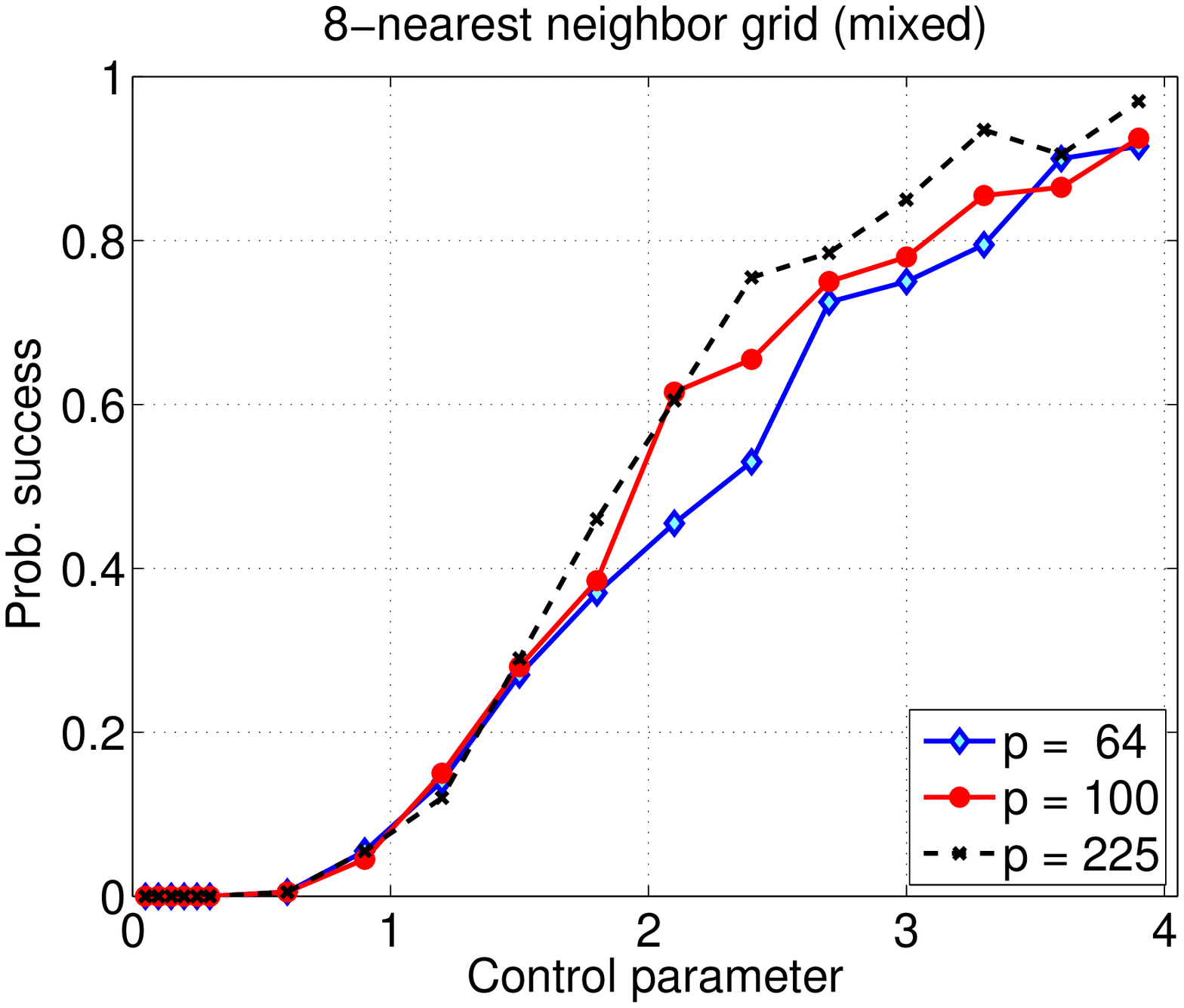} &
\widgraph{0.45\textwidth}{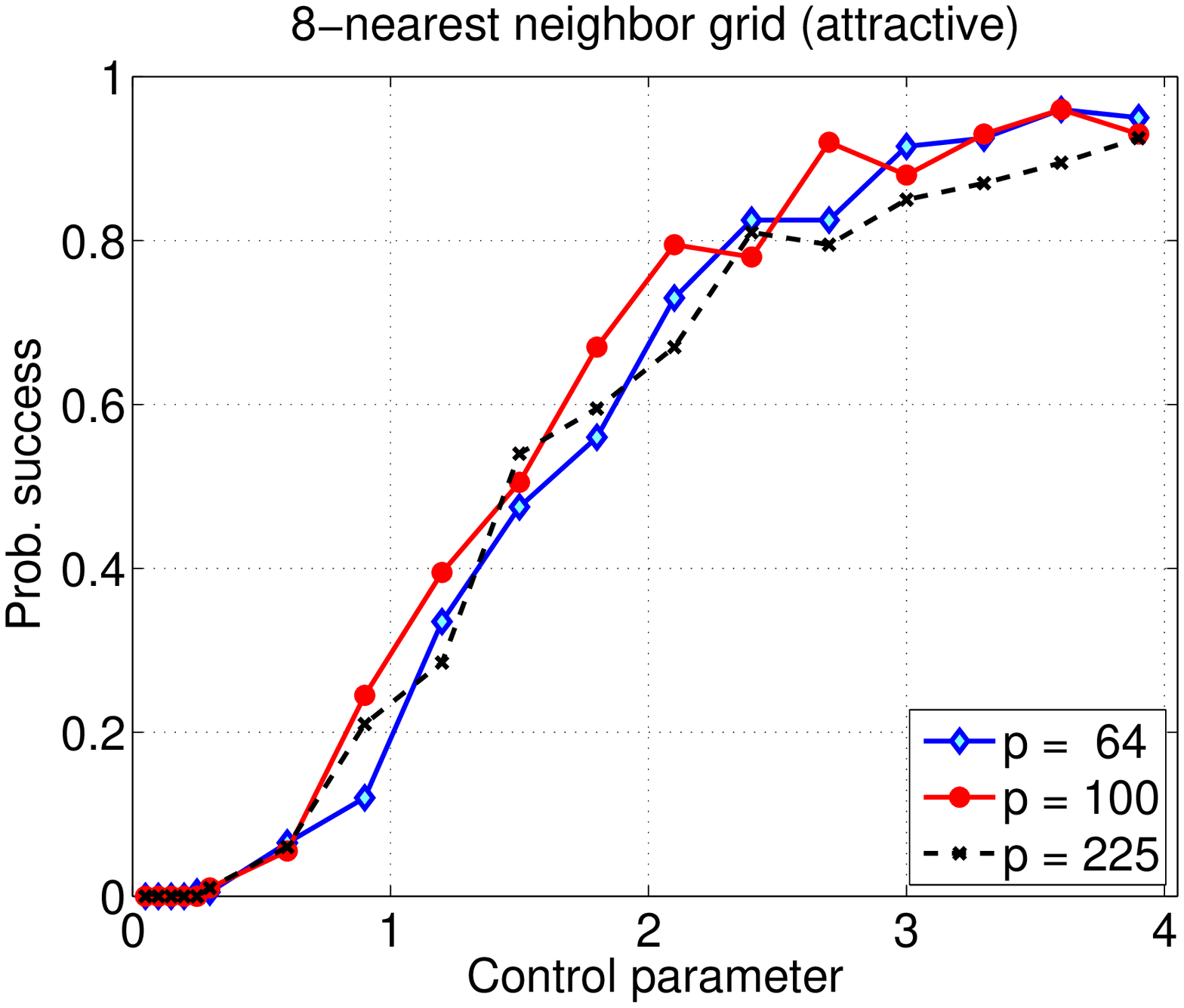} \\
(a) & (b)
\end{tabular}
\caption{Plots of success probability $\mprob[\hatNsign(\svert) =
\N(\svert), \forall r]$ versus the control parameter \mbox{$\contpar(\numobs,
\pdim, \degmax) = \numobs/[10 \degmax \log(\pdim)]$} for Ising models
on 2-D grids with eight nearest-neighbor interactions ($\degmax = 8$).
(a) Randomly chosen mixed sign couplings $\eparam^*_{st} = \pm 0.25$.
(b) All positive couplings $\eparam^*_{st} = 0.25$. }
\label{FigIsingEight}
\end{center}
\end{figure}

Figure~\ref{FigIsingFour} shows results for the $4$-nearest-neighbor
grid model, illustrated in Figure~\ref{FigGraphs}(a), for three different
graph sizes $\pdim \in \{64, 100, 225\}$, with mixed couplings (panel
(a)) and attractive couplings (panel (b)).  Each curve corresponds to
a given problem size, and corresponds to the success probability
versus the control parameter $\contpar$.  Each point corresponds to
the average of $N = 200$ trials.  Notice how despite the very
different regimes of $(\numobs, \pdim)$ that underlie each curve, the
different curves all line up with one another quite well.  This fact
shows that for a fixed degree graph (in this case $\deg = 4$), the
ratio $\numobs/\log(\pdim)$ controls the success/failure of our model
selection procedure, consistent with the prediction of
Theorem~\ref{ThmMain}.
Figure~\ref{FigIsingEight} shows analogous results for the
$8$-nearest-neighbor lattice model ($\degmax = 8$), for the same range
of problem size $\pdim \in \{64, 100, 225\}$, as well as both mixed
and attractive couplings.  Notice how once again the curves for
different problem sizes are all well-aligned, consistent with the
prediction of Theorem~\ref{ThmMain}.

\begin{figure}[ht]
\begin{center}
\begin{tabular}{cc}
\widgraph{0.45\textwidth}{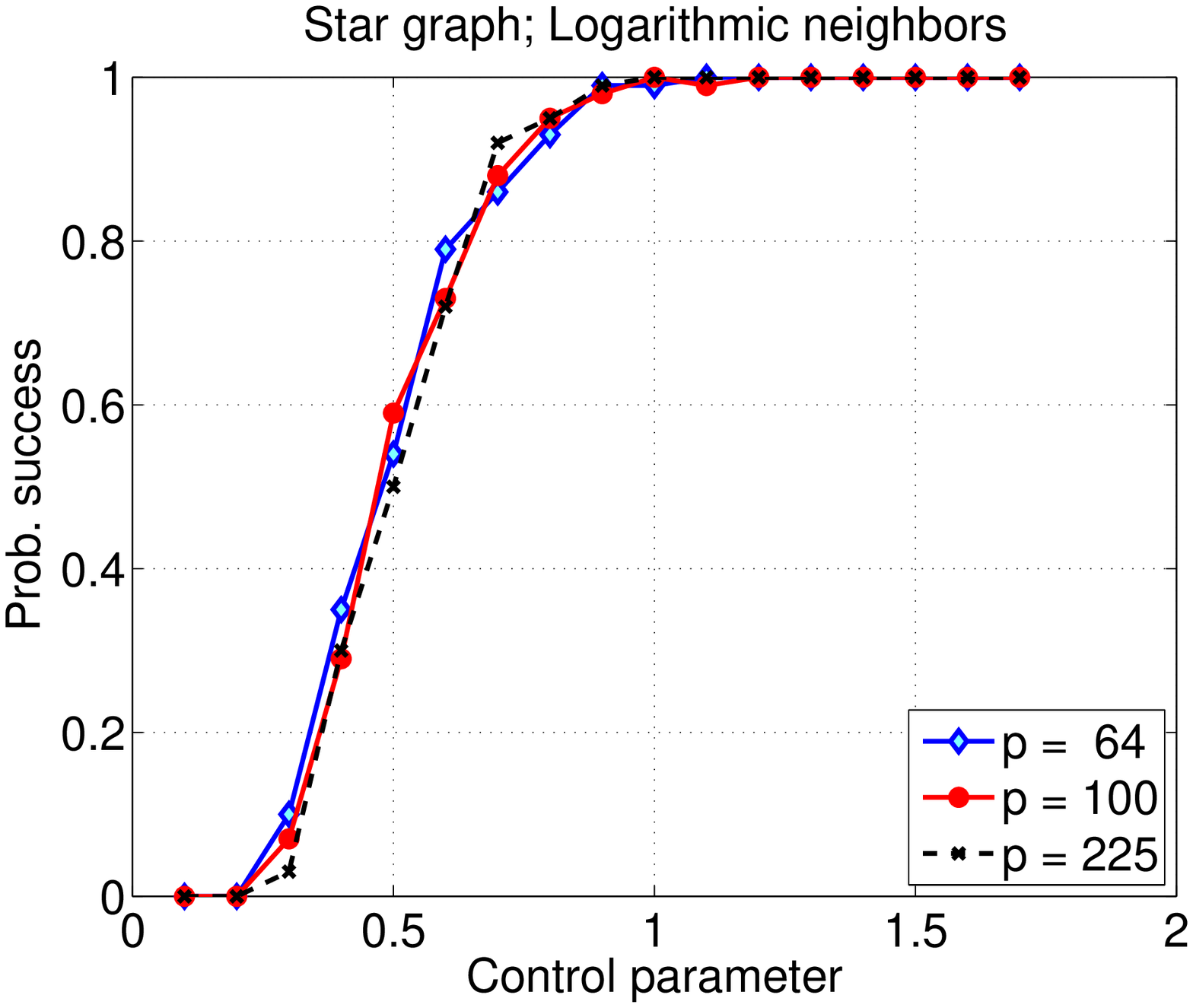}  &
\widgraph{0.45\textwidth}{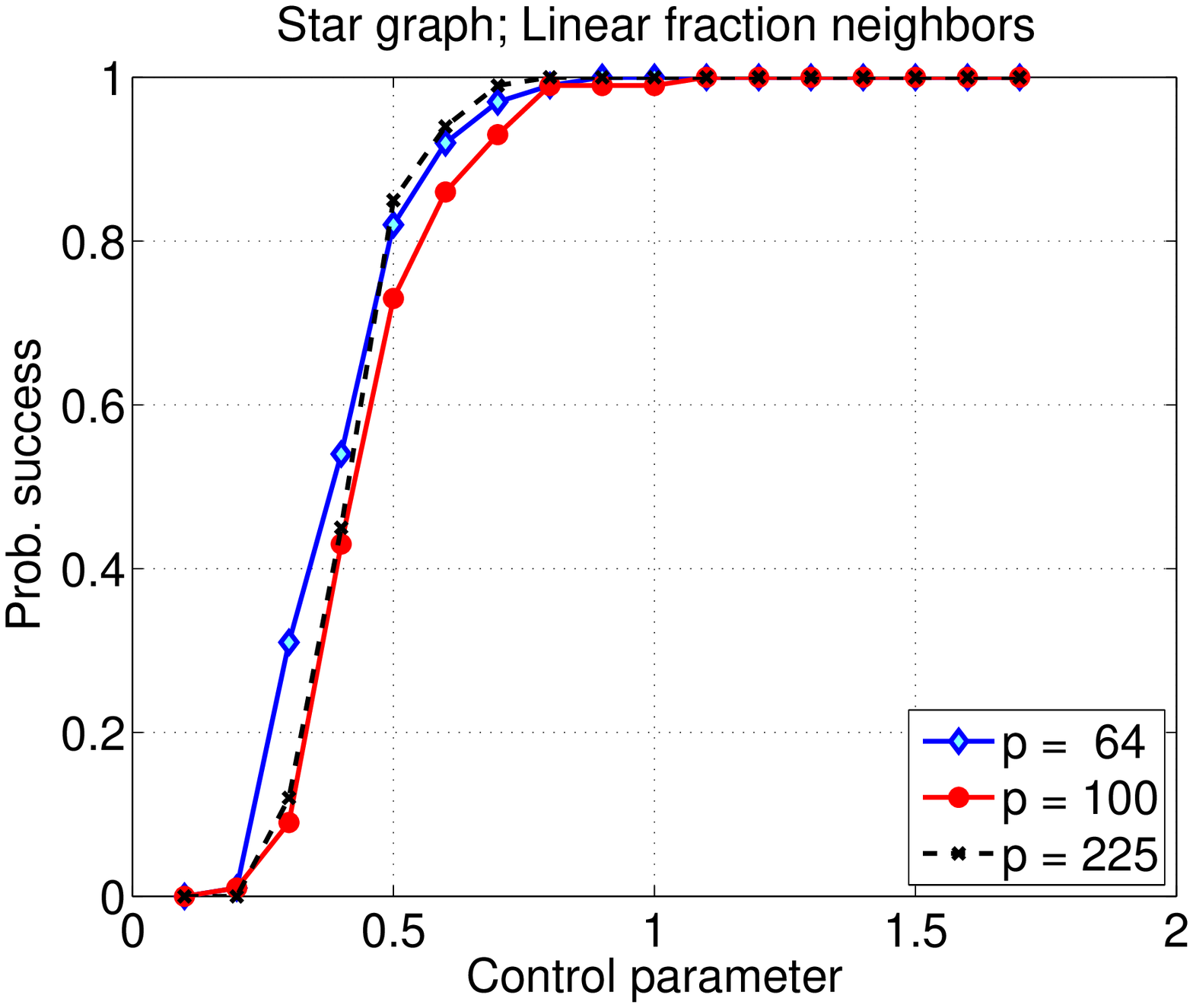}  \\
(a) & (b)
\end{tabular}
\caption{Plots of success probability $\mprob[\hatNsign(\svert) =
\N(\svert), \forall r]$ versus the control parameter \mbox{$\contpar(\numobs,
\pdim, \degmax) = \numobs/[10 \degmax \log(\pdim)]$} for star-shaped
graphs, in which $\degmax = \Theta(\pdim)$, for attractive couplings.
(a) Logarithmic growth in degrees.  (b) Linear growth in degrees.}
\label{FigStar}
\end{center}
\end{figure}

For our last set of experiments, we investigated the performance
of our method for a class of graphs with unbounded maximum degree
$\degmax$.  In particular, we constructed star-shaped graphs with
$\pdim$ vertices by designating one node as the spoke, and connecting
it to $\degmax < (\pdim -1)$ of its neighbors.  For linear
sparsity, we chose $\degmax = \lceil 0.1 \pdim \rceil$, whereas for
logarithmic sparsity we choose $\degmax = \lceil \log(\pdim) \rceil$.
We again studied a triple of graph sizes $\pdim \in \{64, 100, 225 \}$
and Figure~\ref{FigStar} shows the resulting curves of success
probability versus control parameter $\contpar = \numobs/[10 \degmax
\log(\pdim)]$.  Panels (a) and (b) correspond respectively to the
cases of logarithmic and linear degrees. As with the bounded degree
models in Figure~\ref{FigIsingFour} and~\ref{FigIsingEight}, these curves
align with one another, showing a transition from failure to success
with probability one.

\section{Conclusion}
\label{SecDiscussion}

We have shown that a technique based on $\ell_1$-regularized logistic
regression can be used to perform consistent model selection in
discrete graphical models, with polynomial computational complexity
and sample complexity logarithmic in the graph size.  Our analysis
applies to the high-dimensional setting, in which both the number of
nodes $\pdim$ and maximum neighborhood sizes $\degmax$ are allowed to
grow as a function of the number of observations $\numobs$.  There are
a number of possible directions for future work.  For bounded degree
graphs, our results show that the structure can be recovered with high
probability once $\numobs/\log(\pdim)$ is sufficiently large.  Up to
constant factors, this result matches known information-theoretic
lower bounds~\citep{Bresler08,SanWai08}. On the other hand, our
experimental results on graphs with growing degrees (star-shaped
graphs) are consistent with the conjecture that the logistic
regression procedure exhibits a threshold at a sample size $\numobs =
\Theta(\degmax \log \pdim)$, at least for problems where the minimum
value $\thetamin$ stays bounded away from zero. It would be
interesting to provide a sharp threshold result for this problem,
to parallel the known thresholds for $\ell_1$-regularized linear
regression, or the Lasso (see~\cite{Wainwright06a_aller}). 
Finally, the ideas described here, while specialized in this paper to the
pairwise binary case, are more broadly applicable to discrete
graphical models with a higher number of states; this is an
interesting direction for future research.

\ifthenelse{\equal{\doctype}{TECH}}
{
\subsection*{Acknowledgements} Research supported in part by NSF grants IIS-0427206 and CCF-0625879
(PR and JL), NSF grants DMS-0605165 and CCF-0545862 (PR and MJW), and
a Siebel Scholarship~(PR).
}{}
\ifthenelse{\equal{\doctype}{JMLR}}
{\acks{Research supported in part by
NSF grants IIS-0427206 and CCF-0625879 (PR and JL), NSF grants
DMS-0605165 and CCF-0545862 (PR and MJW), and a Siebel Scholarship
(PR).}  
}{}
%%%%%%%%%%%%%%%%%%%%%%%%%%%%%%%%%%%%%%%%%%%%%%%%%%%%%%%%%%%%%%%%%%%%%%%%%%%%%%

\appendix

\section{Proofs for Section~4.1}
\label{AppTechnical}
In this section, we provide proofs of Lemmas~\ref{LemBernOne},
Lemma~\ref{Leml2cons} and Lemma~\ref{LemTayRem}, previously stated in
Section~\ref{SecTechnical}.

\subsection{Proof of Lemma~\ref{LemBernOne}}
\label{AppBernOne}

 Note that any entry of $\obsnoise$ has the form $\obsnoise_u =
\frac{1}{\numobs} \sum_{i=1}^\numobs \zsam{i}_u$, where for $i=1, 2,
\ldots, \numobs$, the variables
\begin{eqnarray*}
\zsam{i}_u & \defn & \myxsam{i}_{\bk \svert} \left\{\myxsam{i}_\svert
- \mprob_{\eparams}[x_\svert = 1 \, \mid \, \myxsam{i}_{\bk \svert}] +
\mprob_{\eparams}[x_\svert = -1 \, \mid \, \myxsam{i}_{\bk \svert}]
\right\}
\end{eqnarray*}
are zero-mean under $\mprob_{\eparams}$, i.i.d., and bounded
($|\zsam{i}_u| \leq 2$).  Therefore, by the Azuma-Hoeffding
inequality~\citep{Hoeffding63}, we have, for any $\delta > 0$, $\mprob
\left[|\obsnoise_u| > \delta \right] \leq 2 \exp \left( -
\frac{\numobs \delta^2}{8} \right)$.  Setting $\delta = \frac{\mutinc
\regpar_\numobs}{4\,(2-\mutinc)}$, we obtain
\begin{eqnarray*}
\mprob\left[\frac{2-\mutinc}{\regpar_\numobs} |\obsnoise_u| >
\frac{\mutinc}{4}\right] & \leq & 2 \exp \left( - \unicon \numobs
\regpar^2_\numobs \right)
\end{eqnarray*}
for some constant $\unicon$.  Finally, applying a union bound 
over the indices $u$ of $\obsnoise$ yields
\begin{eqnarray*}
\mprob\left[\frac{2-\mutinc}{\regpar_\numobs} \|\obsnoise\|_\infty >
\frac{\mutinc}{4}\right] & \leq & 2 \exp \left( - \unicon \numobs
\regpar^2_\numobs + \log(\pdim)\right),
\end{eqnarray*}
as claimed.

\subsection{Proof of Lemma~\ref{Leml2cons}}
\label{Appl2cons}
Following a method of~\cite{Rothman08}, we define the function
$G:\real^{\degmax} \rightarrow \real$ by
\begin{align}
G(\pert_\Nset) & \mydefn  \neglog(\tparam_{\Nset} + \pert_\Nset;
\Data) - \neglog(\tparam_{\Nset}; \Data) +
\regpar_\numobs\left(\|\tparam_{\Nset} + \pert_\Nset\| -
\|\tparam_{\Nset}\|\right).
\end{align}
It can be seen from equation~\eqref{thetaSCon} that $\estim{\pert} =
\hat{\eparam}_{S} - \tparam_{S}$ minimizes $G$.  Moreover $G(0) = 0$ by 
construction; therefore, we must have $G(\estim{\pert}) \leq 0$.
Note also that $G$ is convex. Suppose that we show that for some
radius $B > 0$, and for $\pert \in \real^\degmax$ with $\|\pert\|_2 =
B$, we have $G(\pert) > 0$.  We then claim that $\|\estim{\pert}\|_2
\leq B$.  Indeed, if $\estim{\pert}$ lay outside the ball of radius
$B$, then the convex combination $t \estim{\pert} + (1-t) (0)$ would
lie on the boundary of the ball, for an appropriately chosen $t \in
(0,1)$.  By convexity,
\begin{eqnarray*}
G\left(t \estim{\pert} + (1-t) (0)
\right) \; \leq \; tG(\estim{\pert}) + (1-t) G(0)  \, \leq \, 0,
\end{eqnarray*}
contradicting the assumed strict positivity of $G$ on the boundary.

It thus suffices to establish strict positivity of $G$ on the boundary
of the ball with radius \mbox{$B = \Mrad \regpar_\numobs
\sqrt{\degmax}$,} where $\Mrad > 0$ is a parameter to be chosen later
in the proof.  Let $\pert \in \real^\degmax$ be an arbitrary vector
with $\|\pert\|_2 = B$.  Recalling the notation $\Wnoise = \nabla
\neglog(\eparams; \Data)$, by a Taylor series expansion of the log
likelihood component of $G$, we have
\begin{eqnarray}
\label{EqnGtaylor}
G(\pert) & = & \Wnoise_{\Nset}^T \pert + \pert^T \left[ \Hess
\neglog(\tparam_{\Nset} + \alpha \pert)\right] \pert +
\regpar_\numobs\left(\|\tparam_{\Nset} + \pert_{S}\| -
\|\tparam_{S}\|\right),
\end{eqnarray}
for some $\alpha \in [0,1]$.
For the first term, we have the bound
\begin{eqnarray}
\label{EqnBoundOne}
|\Wnoise_\Nset^T \pert| & \leq & \|\Wnoise_\Nset\|_\infty \|\pert\|_1
 \; \leq \; \|\Wnoise_\Nset\|_\infty \sqrt{\degmax} \, \|\pert\|_2 \;
 \leq \; \left(\regpar_\numobs \sqrt{\degmax} \right)^2
 \frac{\Mrad}{4},
\end{eqnarray}
since $\|\Wnoise_\Nset\|_\infty \leq \frac{\regpar_\numobs}{4}$ with
probability converging to one from Lemma~\ref{LemBernOne}.

Applying the triangle inequality to the last term in the
expansion~\eqref{EqnGtaylor} yields
\begin{align}
\label{EqnBoundTwo}
\regpar_\numobs \|\tparam_{\Nset} + \pert_{\Nset}\|_{1} -
\|\tparam_{\Nset}\|_{1} & \geq  - \regpar_\numobs \|\pert_\Nset\|_1
\, \geq \; -\regpar_\numobs \sqrt{\degmax} \|\pert_\Nset\|_2 \; = \;
-\Mrad \left(\sqrt{\degmax} \regpar_\numobs \right)^2.
\end{align}

Finally, turning to the middle Hessian term, we have
\begin{align*}
q^* & \defn\; \myeigmin(\Hess\neglog(\tparam_{S} + \alpha \pert; \Data))\\
 &\geq\;  \min_{\alpha \in [0,1]} \myeigmin(\Hess\neglog(\tparam_{S} +
\alpha \pert_{S}; \Data))\\
 &=\; \min_{\alpha \in [0,1]}
\myeigmin\left[\frac{1}{n}\sum_{i=1}^{n}\varfun{\xi}{\tparam_S +
\alpha \pert_S} \xi_{S} (\xi_{S})^T\right]
\end{align*}
By a Taylor series expansion of $\varfun{\xi}{\cdot}$, we have
\begin{align*}
q^* & \geq\; \myeigmin\left[\frac{1}{n}\sum_{i=1}^{n}\varfun{\xi}{\tparam_S}\xi_{S}
(\xi_{S})^T\right]\\
&\qquad\quad - \max_{\alpha \in [0,1]}
\matsnorm{\frac{1}{n}\sum_{i=1}^{n}\eta'(\xi; \tparam_\Sset + \alpha
\pert_\Sset)(\pert_{S}^{T}\xi_{S})\xi_{S} (\xi_{S})^T}{2} \\
& =\;  \myeigmin(\Qstar_{SS}) - \max_{\alpha \in [0,1]}
\matsnorm{\frac{1}{n}\sum_{i=1}^{n}\eta'(\xi; \tparam_\Sset + \alpha
\pert_\Sset)(\pert_{S}^{T}\xi_{S})\xi_{S} (\xi_{S})^T}{2} \\
& \geq  \Cmin - \max_{\alpha \in [0,1]}
\matsnorm{\frac{1}{n}\sum_{i=1}^{n}\eta'(\xi; \tparam_\Sset + \alpha
\pert_\Sset)(\pert_{S}^{T}\xi_{S})\xi_{S} (\xi_{S})^T}{2} \\
\end{align*}
It remains to control the final spectral norm.  For any fixed $\alpha
\in [0,1]$ and $y \in \real^\degmax$ with $\|y \|_2 = 1$, we have
\begin{eqnarray*}
\lefteqn{y^T \left \{\frac{1}{n}\sum_{i=1}^{n}\eta'(\xi; \tparam_\Sset + \alpha
\pert_\Sset)(\pert_{S}^{T}\xi_{S})\xi_{S} (\xi_{S})^T \right \} y \;
=\;} \\
&& \frac{1}{n}\sum_{i=1}^{n} \eta'(\xi; \tparam_\Sset + \alpha
\pert_\Sset)(\pert_{S}^{T}\xi_{S}) \; \left[\xi_{S})^T y \right]^2 \\
& \leq & \frac{1}{n}\sum_{i=1}^{n} \left |\eta'(\xi; \tparam_\Sset +
\alpha \pert_\Sset)(\pert_{S}^{T}\xi_{S}) \right| \; \left[\xi_{S})^T
y \right]^2 
\end{eqnarray*}
Now note that $|\eta'(\xi; \tparam_\Sset + \alpha \pert_\Sset) | \leq
1$, and $|\pert_{S}^{T}\xi_{S}| \le \sqrt{\degmax } \|\pert_{S}\|_2 =
\Mrad \regpar_\numobs \degmax$.  Moreover, we have
$\|\frac{1}{\numobs} \sum_{i=1}^{n} \left(\xi_{S})^T y\right) \leq
\|\frac{1}{\numobs} \sum_{i=1}^{n}\xi_{S}(\xi_{S})^T\|_2 \; \leq \;
\Dmax$ by assumption.  Combining these pieces, we obtain
\begin{eqnarray*}
\max_{\alpha \in [0,1]} \matsnorm{\frac{1}{n}\sum_{i=1}^{n}\eta'(\xi;
\tparam_\Sset + \alpha \pert_\Sset)(\pert_{S}^{T}\xi_{S})\xi_{S}
(\xi_{S})^T}{2} & \leq & \Dmax \Mrad \regpar_\numobs \degmax \\
& \leq &  \frac{\Cmin}{2},
\end{eqnarray*}
where the last inequality follows as long as $\regpar_\numobs \degmax
\leq \frac{\Cmin}{2 \Dmax \Mrad}$.  We have thus shown that
\begin{eqnarray}
\label{EqnBoundThree}
q^* \defn \myeigmin(\Hess\neglog(\tparam_{S} + \alpha \pert; \Data)) &
\geq & \frac{\Cmin}{2}
\end{eqnarray}
with probability converging to one, as long as $\regpar_\numobs \degmax$ is
sufficiently small.

Finally, combining the bounds~\eqref{EqnBoundOne},
~\eqref{EqnBoundTwo}, and~\eqref{EqnBoundThree} in the
expression~\eqref{EqnGtaylor}, we conclude that
\begin{eqnarray*}
G(\pert_\Sset) & \geq & \left(\regpar_\numobs \sqrt{\degmax}\right)^2
\, \left \{ -\frac{1}{4} \Mrad + \frac{\Cmin}{2} \Mrad^2 - \Mrad \right \}.
\end{eqnarray*}
This expression is strictly positive for $\Mrad = 5/\Cmin$.  Moreover,
for this choice of $\Mrad$, we have that $\regpar_\numobs \degmax$ must be
upper bounded by $\frac{\Cmin}{2 \Dmax \Mrad} = \frac{\Cmin^2}{10\Dmax}$, as 
assumed in the lemma statement. 

\subsection{Proof of Lemma~\ref{LemTayRem}}
\label{AppTayRem}

We first show that the remainder term $\Rem^\numobs$ satisfies the
bound $\|\Rem^\numobs\|_\infty \leq \Dmax \|\eparamhat_\Sset -
\eparam^*_\Sset\|_2^2$.  Then the result of
Lemma~\ref{Leml2cons}---namely, that \mbox{$\|\eparamhat_\Sset -
\eparam^*_\Sset\|_2 = \mathcal{O}_p(\regpar_\numobs
\sqrt{\degmax})$}---can be used to conclude that\\
\mbox{$\frac{\|\Rem^\numobs\|_\infty}{\regpar_\numobs} =
\mathcal{O}_p(\regpar_\numobs \degmax)$,} which suffices to guarantee
the claim of Lemma~\ref{LemTayRem}.

\newcommand{\myfunc}{\ensuremath{g}} 
\newcommand{\eparambarbar}{\ensuremath{\bar{\bar{\eparam}}}}

Focusing on element $\Rem^\numobs_j$ for some
index $j \in \{1, \ldots, \pdim \}$, we have
\begin{eqnarray*}
\Rem^\numobs_j & = & \left[\nabla^2 \neglog(\eparambar^{(j)}; x) - \nabla^2
\neglog(\eparams; x) \right]_j^T \, [\eparamhat - \eparam^*] \\
& = & \frac{1}{\numobs} \sum_{i=1}^\numobs
\left[\varfun{\xi}{\eparambar^{(j)}} - \varfun{\xi}{\eparams} \right] \left[\xi
(\xi)^T\right]_j^T [\eparamhat - \eparam^*].
\end{eqnarray*}
for some point $\eparambar^{(j)} = t_j\eparamhat + (1-t_j)\eparams$.
Setting $\myfunc(t) = \frac{4
\exp(2t)}{[1+ \exp(2t)]^2}$, note that $\varfun{x}{\eparam} = \myfunc(
x_\svert \sum_{t \in \vertex \bk \svert} \eparam_{\svert t} x_t )$.  
By the chain rule and another application of the 
mean value theorem, we then have
\begin{eqnarray*}
\Rem^\numobs_j & = & \frac{1}{\numobs} \sum_{i=1}^\numobs
\myfunc'\left(\eparambarbar^{(j)T} \xi\right) (\xi)^T [\eparambar^{(j)} - \eparam^*] \left \{
\xi_j (\xi)^T [\eparamhat - \eparam^*] \right \} \\
   & = & \frac{1}{\numobs} \sum_{i=1}^\numobs \left
\{\myfunc'\left(\eparambarbar^{(j)T} \xi\right) \xi_j \right \} \left \{
\eparambar^{(j)}-\eparams]^T \xi (\xi)^T
[\eparamhat - \eparam^*] \right \}
\end{eqnarray*}
where $\eparambarbar^{(j)}$ is another point on the line joining $\eparamhat$ and $\eparams$.
Setting $a_i \defn \{\myfunc'\left(\eparambarbar^{(j)T} \xi\right) \xi_j \}$ and $b_i
\defn \{ [\eparambar^{(j)}-\eparam^*]^T \xi (\xi)^T [\eparamhat - \eparam^*] \}$, we have
\begin{eqnarray*}
|\Rem^\numobs_j| & = & \frac{1}{\numobs} \left|\sum_{i=1}^\numobs a_i
b_i \right| \; \leq \; \frac{1}{\numobs} \|a\|_\infty \|b\|_1.
\end{eqnarray*}
A calculation shows that $\|a\|_\infty \leq 1$, and
\begin{eqnarray*}
\frac{1}{\numobs} \|b\|_1 & = &  t_j [\eparamhat - \eparam^*]^T \left
\{\frac{1}{\numobs} \sum_{i=1}^n \xi (\xi)^T \right \} [\eparamhat -
\eparam^*] \\
& = & t_j [\eparamhat_\Sset - \eparam^*_\Sset]^T \left \{\frac{1}{\numobs}
\sum_{i=1}^n \xi_\Sset (\xi_\Sset)^T \right \} [\eparamhat_\Sset -
\eparam^*_\Sset] \\
& \leq & \Dmax \|\eparamhat_\Sset - \eparam^*_\Sset\|_2^2,
\end{eqnarray*}
where the second line uses the fact that $\eparamhat_\Sbar =
\eparam^*_\Sbar = 0$. This concludes the proof.

\section{Proof of Lemma~7}
\label{AppTechBounds}

Recall from the discussion leading up to the bound~\eqref{EqnAzuma}
that element $(j,k)$ of the matrix difference $\Qobs - \Qstar$,
denoted by $\myZ_{jk}$, satisfies a sharp tail bound.  By definition of
the $\ell_\infty$-matrix norm, we have
\begin{eqnarray*}
\mprob[\matnorm{\Qobs_{\Sbar \Sset} - \Qstar_{\Sbar
\Sset}}{\infty}{\infty} \geq \delta] & = & \mprob[\max_{j \in \Sbar} \sum_{k
\in \Sset} |\myZ_{jk}|  \geq \delta] \\
& \leq & (\pdim-\degmax) \, \mprob[\sum_{k\in \Sset} |\myZ_{jk}| \geq
\delta],
\end{eqnarray*}
where the final inequality uses a union bound, and the fact that
$|\Sbar| \leq \pdim - \degmax$.  Via another union bound over the
row elements, we have
\begin{eqnarray*}
\mprob[\matnorm{\Qobs_{\Sbar \Sset} - \Qstar_{\Sbar
\Sset}}{\infty}{\infty} \geq \delta] & \leq (\pdim - \degmax) \;
\degmax \; \mprob \left[|\myZ_{jk}| \geq \delta/\degmax \right],
\end{eqnarray*}
from which the claim~\eqref{EqnTechBoundA} follows by setting
$\epsilon = \delta/\degmax$ in the Hoeffding bound~\eqref{EqnAzuma}.
The proof of bound~\eqref{EqnTechBoundB} is analogous, with the
pre-factor $(\pdim - \degmax)$ replaced by $\degmax$.

To prove the last claim~\eqref{EqnTechBoundC}, we write
\begin{eqnarray*}
\matnorm{(\Qobs_{\Sset \Sset})^{-1} - (\Qstar_{\Sset
\Sset})^{-1}}{\infty}{\infty} & = & \matsnorm{ (\Qstar_{\Sset
\Sset})^{-1} \left[\Qstar_{\Sset \Sset} - \Qobs_{\Sset \Sset} \right]
(\Qobs_{\Sset \Sset})^{-1} }{\infty} \\
& \leq & \sqrt{\degmax} \; \matsnorm{ (\Qstar_{\Sset
\Sset})^{-1} \left[\Qstar_{\Sset \Sset} - \Qobs_{\Sset \Sset} \right]
(\Qobs_{\Sset \Sset})^{-1} }{2} \\
& \leq & \sqrt{\degmax} \; \matsnorm{ (\Qstar_{\Sset \Sset})^{-1}}{2}
\matsnorm{\Qstar_{\Sset \Sset} - \Qobs_{\Sset \Sset}}{2}
\matsnorm{(\Qobs_{\Sset \Sset})^{-1}}{2} \\
& \leq & \frac{\sqrt{\degmax}}{\Cmin} \matsnorm{\Qstar_{\Sset \Sset} -
\Qobs_{\Sset \Sset}}{2} \matsnorm{(\Qobs_{\Sset \Sset})^{-1}}{2}.
\end{eqnarray*}
From the proof of Lemma~\ref{LemConcentrateEig}, in particular
equation~\eqref{EqnUseful}, we have
\begin{eqnarray*}
\mprob \left[\matsnorm{(\Qobs_{\Sset \Sset})^{-1}}{2} \geq
\frac{2}{\Cmin} \right] & \leq 2 \exp \left(-\unicon\frac{\delta^2
\numobs}{ \degmax^2} + B\log(\degmax) \right)
\end{eqnarray*}
for a constants $B$.  Moreover, from equation~\eqref{EqnUseful},
we have
\begin{eqnarray*}
\mprob[\matsnorm{\Qobs_{\Sset \Sset} -\Qmat_{\Sset \Sset}}{2} \geq
  \delta/\sqrt{\degmax} ] & \leq & 2 \exp \left(-\unicon \frac{\delta^2
  \numobs}{\degmax^3} + 2 \log(\degmax)\right),
\end{eqnarray*}
so that the bound~\eqref{EqnTechBoundC} follows.

%%%%%%%%%%%%%%%%%%%%%%%%%%%%%%%%%%%%%%%%%%%%%%%%%%%%%%%%%%%%%%%%%%%%%%%%%%%
%%% BIBLIOGRAPHY
%\vskip 0.2in \bibliography{mjwain_super}

\end{document}